\documentclass[final]{siamltex}

\usepackage[usenames]{color}
\usepackage[color]{showkeys}
\definecolor{refkey}{gray}{0.4}
\definecolor{labelkey}{gray}{0.3}

\usepackage[left=1in,top=1in,right=1in,bottom=1in,dvips,letterpaper]{geometry}
\usepackage{setspace,url}
\usepackage[numbers,square,sort]{natbib}
\usepackage[leqno]{amsmath}
\usepackage{amsxtra, amsfonts,amscd,  amssymb, graphicx,subfigure}

\usepackage[algo2e, vlined,ruled,linesnumbered]{algorithm2e}
\numberwithin{equation}{section}

\newcommand{\be}{\begin{equation}}
\newcommand{\ee}{\end{equation}}
\newcommand{\ba}{\begin{array}}
\newcommand{\ea}{\end{array}}
\newcommand{\bea}{\begin{eqnarray}}
\newcommand{\eea}{\end{eqnarray}}
\newcommand{\beaa}{\begin{eqnarray*}}
\newcommand{\eeaa}{\end{eqnarray*}}

\newcommand{\half}{\frac{1}{2}}
\newcommand{\br}{\mathbb{R}}

\newcommand{\mtn}{{m \times n}}

\usepackage[usenames]{color}
\usepackage[normalem]{ulem}

\newcommand{\LCal}{\mathcal{L}}
\newcommand{\ACal}{\mathcal{A}}

\newcommand{\CCal}{\mathcal{C}}

\newcommand{\PCal}{\mathcal{P}}

\newcommand{\etal}{{et al. }}
\newcommand{\Diag}{\mbox{Diag}}
\newcommand{\oneotwomu}{\frac{1}{2\mu}}
\newcommand{\oneomu}{\frac{1}{\mu}}
\newcommand{\st}{\mbox{ s.t. }}

\newtheorem{remark}[theorem]{Remark}

\begin{document}

\title{Fast Alternating Linearization Methods for Minimizing the Sum of Two Convex Functions}
\date{}
\author{Donald Goldfarb\footnotemark[1] \and Shiqian Ma\footnotemark[1] \and Katya Scheinberg\footnotemark[2]}
\renewcommand{\thefootnote}{\fnsymbol{footnote}}
\footnotetext[1]{Department of Industrial Engineering and Operations
Research, Columbia University, New York, 10027, USA. \quad Email:
\{goldfarb, sm2756\}@columbia.edu. Research supported in
part by NSF Grants DMS 06-06712 and DMS 10-16571, ONR Grant
N00014-08-1-1118 and DOE Grant DE-FG02-08ER25856. }
\footnotetext[2]{Department of Industrial and Systems Engineering, Lehigh University, Bethlehem, PA 18015-1582, USA. \quad Email: katyas@lehigh.edu}

\renewcommand{\thefootnote}{\arabic{footnote}}

\maketitle \centerline{October 11, 2010}
\begin{abstract}
We present in this paper first-order alternating linearization algorithms based on an alternating direction augmented Lagrangian approach for minimizing the sum of two convex functions. Our basic methods require at most $O(1/\epsilon)$ iterations to obtain an $\epsilon$-optimal solution, while our accelerated (i.e., fast) versions of them require at most $O(1/\sqrt{\epsilon})$ iterations, with little change in the computational effort required at
each iteration. For both types of methods, we present one algorithm that requires both functions to be smooth with Lipschitz continuous gradients and one algorithm that  needs only one of the  functions to be so.
Algorithms in this paper are Gauss-Seidel type methods, in contrast to the ones proposed by Goldfarb and Ma in
\citep{Goldfarb-Ma-Ksplit} where the algorithms are Jacobi type methods. Numerical results
are reported to support our theoretical conclusions and demonstrate the practical potential of our algorithms.
\end{abstract}



\section{Introduction}\label{sec:intro} In this paper, we are interested in the following convex optimization problem: \bea\label{prob:P} \min F(x)\equiv
f(x)+g(x) \eea where $f,g:\br^n\rightarrow\br$ are both convex
functions such that the following two problems are easy to solve for any $\tau>0$ and $z\in\br^n$ relative to minimizing $F(x)$:
\bea\label{prob:f-shrink}\min\left\{\tau f(x)+\half\|x-z\|^2 \right\}\eea and \bea\label{prob:g-shrink}\min\left\{\tau g(x)+\half\|x-z\|^2 \right\}.\eea
  In particular, we are specially interested in cases where solving \eqref{prob:f-shrink} (or \eqref{prob:g-shrink}) takes roughly the same effort as  computing  the gradient (or a subgradient) of $f(x)$  (or $g(x)$, respectively).  Problems of this type arise in many applications of practical interest. The following are some interesting examples.

{\bf Example 1. $\ell_1$ minimization in compressed sensing (CS). }
Signal recovery problems in compressed sensing
\citep{Candes-Romberg-Tao-2006,Donoho-06} use the $\ell_1$ norm $\|x\|_1:=\sum_{i=1}^n|x_i|$ as a regularization term to enforce sparsity in the solution $x\in\br^n$ of a linear system $Ax=b$, where $A\in\br^\mtn$ and $b\in\br^m$. This results in the unconstrained problem
\bea\label{prob:CS}\min\left\{\half\|Ax-b\|_2^2+\rho\|x\|_1\right\},\eea where $\rho>0$, which is of the form of \eqref{prob:P} with $f(x)=\half\|Ax-b\|_2^2$ and  $g(x):=\rho\|x\|_1$. In this case, the two problems \eqref{prob:f-shrink} and \eqref{prob:g-shrink} are easy to solve. Specifically, \eqref{prob:f-shrink} reduces to solving a linear system and \eqref{prob:g-shrink} reduces to a vector shrinkage operation which requires $O(n)$ operations (see e.g., \citep{Hale-Yin-Zhang-SIAM-2008}).  Depending on the size  and structure of $A$
solving the system of linear equations required by
\eqref{prob:f-shrink} may be more expensive, less expensive or comparable to computing the gradient $A^\top (Ax-b)$ of $f(x)$. In the application we consider in Section \ref{sec:comparison}  these computations are comparable due to the special structure of $A$.

{\bf Example 2. Nuclear norm minimization (NNM). } The nuclear norm minimization problem, which seeks a low-rank solution of a linear system, can be cast as \bea\label{prob:NNM}\min\left\{\half\|\ACal(X)-b\|_2^2+\rho\|X\|_*\right\},\eea where $\tau>0$, $X\in\br^\mtn$, $\ACal:\br^\mtn\rightarrow\br^p$ is a linear operator, $b\in\br^p$ and the nuclear norm $\|X\|_*$ is defined as the sum of the singular values of the matrix $X$. Problem \eqref{prob:NNM} and a special case of it, the so-called matrix completion problem, have many applications in optimal control, online recommendation systems, computer vision, etc. (see e.g., \citep{Recht-Fazel-Parrilo-2007,Candes-Recht-2008,Candes-Tao-2009,Keshavan-Montanari-Oh-2009-1}). In Problem \eqref{prob:NNM}, if we let $f(X)=\half\|\ACal(X)-b\|_2^2$ and $g(X)=\rho\|X\|_*$, then problem \eqref{prob:f-shrink} reduces to solving a linear system. Problem \eqref{prob:g-shrink} has a closed-form solution that is given by matrix shrinkage operation (see e.g., \citep{Ma-Goldfarb-Chen-2008}).

{\bf Example 3. Robust principal component analysis (RPCA). } The RPCA problem seeks to recover a low-rank matrix $X$ from a corrupted matrix $M$. This problem has many applications in computer vision, image processing and web data ranking (see e.g., \citep{Candes-Li-Ma-Wright-RPCA-2009}), and can be formulated as
\bea\label{prob:RPCA}\min\left\{\|X\|_*+\rho\|Y\|_1: X+Y=M\right\},\eea where $\rho>0$, $M\in\br^\mtn$ and the $\ell_1$ norm $\|Y\|_1:=\sum_{i,j}|Y_{ij}|$. Note that \eqref{prob:RPCA} can be rewritten as \beaa\min\{\|X\|_*+\rho\|M-X\|_1\},\eeaa which is of the form of \eqref{prob:P}. Moreover, the two problems \eqref{prob:f-shrink} and \eqref{prob:g-shrink} corresponding to \eqref{prob:RPCA} have closed-form solutions given respectively by a matrix shrinkage operation and a vector shrinkage operation.  The matrix shrinkage operation requires
a singular value decomposition (SVD) and is  comparable in cost to computing a subgradient of $\|X\|_*$ or the gradient of the smoothed version of this function (see  Section  \ref{sec:application}).

{\bf Example 4. Sparse inverse covariance selection (SICS). } Gaussian graphical models are of great interest in statistical learning. Because conditional independence between different nodes correspond to zero entries in the inverse covariance matrix of the Gaussian distribution, one can learn the structure of the graph by estimating a sparse inverse covariance matrix from sample data by solving the following maximum likelihood problem with an $\ell_1$-regularization term, (see e.g., \citep{Yuan-Lin-2007,Friedman-Hastie-Tibshirani-2007,Wainwright-Ravikumar-Lafferty-2007,Banerjee-ElGhaoui-Aspremont-2008}). \beaa \max\left\{\log\det(X)-\langle\Sigma,X\rangle -\rho\|X\|_1\right\},\eeaa or equivalently, \bea\label{prob:SICS}\min\left\{-\log\det(X)+\langle\Sigma,X\rangle+\rho\|X\|_1\right\},\eea where $\rho>0$ and $\Sigma\in S_+^n$ (the set of symmetric positive semidefinite matrices) is the sample covariance matrix. Note that by defining $f(X):=-\log\det(X)+\langle\Sigma,X\rangle$ and $g(X):=\rho\|X\|_1$, \eqref{prob:SICS} is of the form of \eqref{prob:P}. Moreover, it can be proved that the problem \eqref{prob:f-shrink} has a closed-form solution, which is given by a spectral decomposition - a comparable effort to computing the gradient of $f(X)$, while the solution of problem \eqref{prob:g-shrink} corresponds to a vector shrinkage operation.

Algorithms for solving problem \eqref{prob:P} have been studied
extensively in the literature. For large-scale problems, for which problems \eqref{prob:f-shrink} and \eqref{prob:g-shrink} are relatively easy to solve, the class of alternating direction methods that are based on variable splitting combined with the
augmented Lagrangian method are particularly important. In these methods, one splits
the variable $x$ into two variables, i.e., one introduces a new variable
$y$ and rewrites Problem \eqref{prob:P} as \bea \label{prob:P-x-y}
\min \{f(x)+g(y): x-y=0 \}. \eea Since Problem \eqref{prob:P-x-y} is an
equality constrained problem, the augmented Lagrangian method can be
used to solve it. Given a penalty parameter $1/\mu$, at the
$k$-th iteration, the augmented Lagrangian method minimizes the augmented Lagrangian function
\bea\label{Lag-function}\LCal_\mu(x,y;\lambda):=f(x)+g(y)-\langle
\lambda,x-y\rangle + \frac{1}{2\mu}\|x-y\|^2,\eea with respect to $x$ and $y$, i.e., it solves the subproblem
\bea \label{Aug-Lag-old-x-y} (x^k,y^k) & :=\arg\min\limits_{x,y}
\LCal_\mu(x,y;\lambda^k) \eea and then updates the Lagrange multiplier $\lambda^k$
via: \bea\label{Aug-Lag-old-lambda} \lambda^{k+1} & := \lambda^k -
\frac{1}{\mu}(x^k-y^k).\eea Minimizing $\LCal_\mu(x,y;\lambda)$ with respect to $x$ and $y$ jointly is often not easy. In fact, it certainly is not any easier than solving the original problem \eqref{prob:P}. However, if one minimizes $\LCal_\mu(x,y;\lambda)$ with respect to $x$ and $y$ alternatingly, one needs to solve problems of the form \eqref{prob:f-shrink} and \eqref{prob:g-shrink}, which as we have already discussed, is often easy to do.
Such an alternating direction augmented Lagrangian method (ADAL) for
solving \eqref{prob:P-x-y} is given below as Algorithm
\ref{alg:Aug-Lag-classical}.
\begin{algorithm2e}\caption{Alternating Direction Augmented Lagrangian Method (ADAL)}
    \label{alg:Aug-Lag-classical}
\linesnumberedhidden \dontprintsemicolon Choose $\mu$, $\lambda^0$ and $x^0=y^0$.\;
\For{$k=0,1,\cdots$}{
$x^{k+1} := \arg\min_{x} \LCal_\mu(x,y^k;\lambda^k)$ \;
$y^{k+1} := \arg\min_{y} \LCal_\mu(x^{k+1},y;\lambda^k)$ \;
$\lambda^{k+1} := \lambda^k-\oneomu (x^{k+1}-y^{k+1})$ \;
}
\end{algorithm2e}

The history of alternating direction methods (ADMs) goes back to the 1950s for solving PDEs \citep{Douglas-Rachford-56,Peaceman-Rachford-55} and to the
1970s for solving variational problems associated with PDEs \citep{Gabay-Mercier-1976,Glowinski-LeTallec-89}. ADMs have
also been applied to solving variational inequality problems by
Tseng \citep{Tseng-1991,Tseng-1990} and He \etal
\citep{He-Liao-Han-Yang-2002,He-Yang-Wang-2000}. Recently, with the
emergence of compressive sensing and subsequent great interest in $\ell_1$ minimization
\citep{Candes-Romberg-Tao-2006,Donoho-06}, ADMs have been applied to
$\ell_1$ and total variation regularized problems arising from
signal processing and image processing. The papers of Goldstein and Osher
\citep{Goldstein-Osher-08}, Afonso \etal
\citep{Afonso-BD-Figueiredo-2009} and Yang and Zhang
\citep{Yang-Zhang-2009} are based on the alternating direction
augmented Lagrangian framework (Algorithm
\ref{alg:Aug-Lag-classical}), and demonstrate that ADMs are very
efficient for solving $\ell_1$ and TV regularized problems. The work
of Yuan \citep{Yuan-2009} and Yuan and Yang \citep{Yuan-Yang-2009}
showed that ADMs can also efficiently solve $\ell_1$-regularized problems arising
from statistics and data analysis. More recently, Wen,
Goldfarb and Yin \citep{Wen-Goldfarb-Yin-2009} and Malick \etal \citep{Malick-Povh-Rendl-Wiegele-2009} applied alternating direction augmented
Lagrangian methods to solve semidefinite programming (SDP) problems. The results in \citep{Wen-Goldfarb-Yin-2009} show that these methods greatly
outperform interior point methods on several classes of well-structured SDP
problems. Furthermore, He \etal proposed an alternating direction based contraction method for solving separable linearly constrained convex problems \citep{He-Tao-Xu-Yuan-2009}.

Another important and related class of
algorithms for solving \eqref{prob:P} is based on operator-splitting. The aim of these
algorithms is to find an $x$ such that
\bea\label{prob:sum-2-maximal-monotone-operator}0\in
T_1(x)+T_2(x),\eea where $T_1$ and $T_2$ are maximal monotone
operators. This is a more general problem than \eqref{prob:P} and ADMs for it have been the focus of a substantial amount of research; e.g., see \citep{Combettes-Pesquet-DR-2007,Eckstein-Bertsekas-1992,Combettes-2004,Eckstein-Svaiter-2008,Combettes-Wajs-05,Lions-Mercier-79,Spingarn-1983}. Since the first-order optimality conditions for \eqref{prob:P} are: \bea\label{prob:zero-two-subgrad} 0
\in \partial f(x)+\partial g(x),\eea where $\partial f(x)$ denotes
the subdifferential of $f(x)$ at the point $x$, a solution to Problem \eqref{prob:P} can be obtained by solving
Problem \eqref{prob:sum-2-maximal-monotone-operator}. For example, see \citep{Eckstein-Bertsekas-1992,Lions-Mercier-79,Spingarn-1983} and references therein for more information on this class of algorithms.

While global convergence results for various splitting and alternating
direction algorithms have been established under
appropriate conditions, our interest here is on iteration complexity bounds for such algorithms. By an iteration complexity bound we mean a bound on the number of iterations needed to obtain an $\epsilon$-optimal solution which is defined as follows.
\begin{definition}
$x_\epsilon\in\br^n$ is called an $\epsilon$-optimal solution to \eqref{prob:P} if $F(x_\epsilon)-F(x^*)\leq\epsilon$, where $x^*$ is an optimal solution to \eqref{prob:P}.
\end{definition}

Complexity bounds for first-order methods for solving convex optimization problems have been given by Nesterov and many others.
In \citep{Nesterov-1983,NesterovConvexBook2004}, Nesterov gave first-order
algorithms for solving smooth unconstrained convex minimization problems with an iteration complexity of $O(\sqrt{L/\epsilon})$, where $L$ is the Lipschitz constant of the gradient of the objective function, and showed that this is the best complexity that is obtainable when only first-order information is used. These methods can
be viewed as accelerated gradient methods where a combination of
past iterates are used to compute the next iterate.
Similar techniques were
then applied to nonsmooth problems
\citep{Nesterov-2005,Tseng-2008,Beck-Teboulle-2009,Nesterov-07} and corresponding
optimal complexity results were obtained. The ISTA (Iterative Shrinkage/Thresholding Algorithm) and FISTA (Fast Iterative Shrinkage/Thresholding Algorithm) algorithms proposed by Beck and Teboulle in \citep{Beck-Teboulle-2009} are designed for solving \eqref{prob:P} when one of the functions (say $f(x)$) is smooth and the other is not.
It is proved in \citep{Beck-Teboulle-2009} that the number of iterations required by ISTA and FISTA to get an $\epsilon$-optimal solution to problem \eqref{prob:P}
are respectively $O(L(f)/\epsilon)$ and $O(\sqrt{L(f)/\epsilon})$, under the assumption that $\nabla f(x)$ is Lipschitz continuous with Lipschitz constant $L(f)$, i.e., \beaa \|\nabla f(x) - \nabla f(y) \|_2 \leq L(f) \|x-y\|_2, \quad \forall x,y\in\br^n.\eeaa
ISTA computes a sequence $\{x^k\}$ via the iteration
\bea\label{alg:ISTA}x^{k+1}:=\arg\min_x Q_f(x,x^k),\eea where \begin{equation}\label{def:Qf}
Q_f(u,v) := g(u)+f(v)+\langle \nabla f(v),u-v\rangle+\oneotwomu\|u-v\|^2,
\end{equation}
while FISTA computes $\{x^k\}$ via the iteration
\bea\label{alg:FISTA}\left\{\ba{lll} x^k & := & \arg\min_x Q_f(x,y^k) \\ t_{k+1} & := & \left(1+\sqrt{1+4t_k^2}\right)/2 \\ y^{k+1} & := & x^k + \left(\frac{t_k-1}{t_{k+1}}\right)(x^k-x^{k-1}) \ea\right.\eea
starting with $t_1=1$, $y^1=x^0\in\br^n$ and $k=1$.

Note that ISTA and FISTA treat the functions $f(x)$ and $g(x)$ very differently. At each iteration they both
linearize the function $f(x)$ but never directly minimize it, while they do minimize the function $g(x)$ in conjunction with the linearization of $f(x)$ and a proximal (penalty) term. These two methods have proved to be efficient for solving the CS problem \eqref{prob:CS} (see e.g., \citep{Beck-Teboulle-2009,Hale-Yin-Zhang-SIAM-2008}) and the NNM problem \eqref{prob:NNM} (see e.g., \citep{Ma-Goldfarb-Chen-2008,Toh-Yun-2009}). ISTA and FISTA work well in these areas because $f(x)$ is quadratic and is well approximated by linearization. However, for the RPCA problem \eqref{prob:RPCA} where two complicated functions are involved, ISTA and FISTA do not work well. As we shall show in Sections \ref{sec:ALM} and \ref{sec:FALM}, our ADMs are very effective in solving RPCA problems. For the SICS problem \eqref{prob:SICS}, intermediate iterates $X^k$ may not be positive definite, and hence the gradient of $f(X)=-\log\det(X)+\langle\Sigma,X\rangle$ may not be well defined at $X^k$. Therefore, ISTA and FISTA cannot be used to solve the SICS problem \eqref{prob:SICS}. In \citep{Scheinberg-Ma-Goldfarb-NIPS-2010}, it is shown that SICS problems can be very efficiently solved by our ADM approach.

{\bf Our contribution.}
In this paper, we propose both basic and accelerated (i.e., {\it fast}) versions of first-order alternating linearization methods (ALMs) based on an alternating
direction augmented Lagrangian approach for solving \eqref{prob:P} and analyze their
iteration complexities. Our basic methods require at most $O(L/\epsilon)$ iterations to obtain an $\epsilon$-optimal solution, while our fast methods require at most $O(\sqrt{L/\epsilon})$ iterations with only a very small increase in the computational effort required at
each iteration. Thus, our fast methods are {\em optimal} first-order methods in terms of iteration complexity. For both types of methods, we present an algorithm that requires both functions to be continuously differentiable with Lipschitz constants for the gradients denoted by $L(f)$ and $L(g)$. In this case $L=\max\{L(f), L(g)\}$. We also present for each type of method, an algorithm that only needs one of the functions, say $f(x)$, to be smooth, in which case $L=L(f)$.  These algorithms are related to
the multiple splitting algorithms in a recent paper by Goldfarb and Ma \citep{Goldfarb-Ma-Ksplit}. The algorithms in \citep{Goldfarb-Ma-Ksplit} are Jacobi type methods since they do not use information from the current iteration to solve
succeeding subproblems in that iteration, while the
algorithms proposed in this paper are Gauss-Seidel type methods
since information from the current iteration is used later in
the same iteration.   These algorithms can also be viewed as extensions of the ISTA and FISTA algorithms in  \citep{Beck-Teboulle-2009}. The complexity bounds we obtain for our algorithms are similar to (and as much as a factor of two better that) those in  \citep{Beck-Teboulle-2009}.

 At each iteration, our algorithms alternatively
minimize two different approximations to the original objective function, obtained by keeping one function unchanged and linearizing the other one. Our basic algorithm is similar in many ways to the alternating linearization method proposed by Kiwiel \etal \citep{Kiwiel-Rosa-1999}. In particular, the approximate functions minimized at each step of Algorithm 3.1 in \citep{Kiwiel-Rosa-1999} have the same form as those minimized in our algorithm. However, our basic algorithm differs from the one in \citep{Kiwiel-Rosa-1999} in the way that the proximal terms are chosen, and our accelerated algorithms are very different. Moreover, no complexity bounds have been given for the algorithm in \citep{Kiwiel-Rosa-1999}. To the best of our knowledge, the complexity results in this paper are the first ones that have been given for a Gauss-Seidel type alternating direction method \footnote{After completion of an earlier version of the present paper, which is available on http://arxiv.org/abs/0912.4571, Monteiro and Svaiter \citep{Monteiro-Svaiter-2010a} gave an iteration complexity bound to achieve a desired closeness of the current iterate to the solution for ADMs for solving the more general problem \eqref{prob:sum-2-maximal-monotone-operator}.}. Complexity results for related Jacobi type alternating direction methods are given in \citep{Goldfarb-Ma-Ksplit}.

{\bf Organization.} The rest of this paper is organized as follows. In Sections
\ref{sec:ALM} and \ref{sec:FALM} we propose our alternating linearization methods based on alternating
direction augmented Lagrangian methods and give convergence/complexity bounds for them. We compare the performance of our ALMs to other competing first-order algorithms using an image deblurring problem in Section \ref{sec:comparison}. In Section \ref{sec:application}, we apply our ALMs to solve very large RPCA problems arising from background extraction in surveillance video and matrix completion and report the numerical results. Finally, we make some conclusion in Section \ref{sec:conclude}.

\section{Alternating Linearization Methods}\label{sec:ALM}
In iteration of the ADAL method, Algorithm \ref{alg:Aug-Lag-classical}, the Lagrange
multiplier $\lambda$ is updated just once, immediately after the augmented Lagrangian is minimized with respect to $y$. Since the alternating direction approach is meant to be symmetric with respect to $x$ and $y$, it is natural to also update $\lambda$ after solving the
subproblem with respect to $x$. By doing this, we get a symmetric version of the ADAL method. This algorithm is given below as
Algorithm \ref{alg:Aug-Lag-symmetric}.
\begin{algorithm2e}\caption{Symmetric Alternating Direction Augmented Lagrangian Method (SADAL)}
    \label{alg:Aug-Lag-symmetric}
\linesnumberedhidden \dontprintsemicolon Choose $\mu$, $\lambda^0$ and $x^0=y^0$.\;
\For{$k=0,1,\cdots$}{
$x^{k+1} := \arg\min_x \LCal_\mu(x,y^k;\lambda^k)$ \;
$\lambda^{k+\half} := \lambda^k - \oneomu (x^{k+1}-y^k)$ \;
$y^{k+1} := \arg\min_y \LCal_\mu(x^{k+1},y;\lambda^{k+\half})$ \;
$\lambda^{k+1} := \lambda^{k+\half}-\oneomu (x^{k+1}-y^{k+1})$ \;
}
\end{algorithm2e}

This ADAL variant is described and analyzed in \citep{Glowinski-LeTallec-89}. Moreover, it is shown in \citep{Glowinski-LeTallec-89} that Algorithms \ref{alg:Aug-Lag-classical} and \ref{alg:Aug-Lag-symmetric} are equivalent to the Douglas-Rachford \citep{Douglas-Rachford-56} and Peaceman-Rachford \citep{Peaceman-Rachford-55} methods, respectively applied to the optimality condition \eqref{prob:zero-two-subgrad} for problem \eqref{prob:P}.
If we assume that both $f(x)$ and $g(x)$ are differentiable, it follows from the first order optimality conditions for the two subproblems in lines 3 and 5 of Algorithm \ref{alg:Aug-Lag-symmetric} that
\bea\label{gradg,gradf,lambda} \lambda^{k+\half}  =\nabla f(x^{k+1}) \quad\mbox{ and }\quad \lambda^{k+1} = -\nabla g(y^{k+1}).\eea
Substituting \eqref{gradg,gradf,lambda} into Algorithm \ref{alg:Aug-Lag-symmetric}, we get the following alternating linearization method (ALM) which is equivalent to the SADAL method (Algorithm \ref{alg:Aug-Lag-symmetric}) when both $f$ and $g$ are differentiable.
\begin{algorithm2e}\caption{Alternating Linearization Method (ALM)} \label{alg:ALM}
\linesnumberedhidden \dontprintsemicolon Choose $\mu$ and $x^0=y^0$.\;
\For{$k=0,1,\cdots$}{
$x^{k+1} := \arg\min_x Q_g(x,y^k)$ \;
$y^{k+1} := \arg\min_y Q_f(y,x^{k+1})$ \;
}
\end{algorithm2e}
In Algorithm \ref{alg:ALM}, $Q_f(u,v)$ is defined by \eqref{def:Qf} and
\bea\label{def:Qg}Q_g(u,v):=f(u)+g(v)+\langle\nabla g(v),u-v\rangle+\oneotwomu\|u-v\|_2^2. \eea

In Algorithm \ref{alg:ALM}, we alternatively replace the functions $g$ and $f$ by their linearizations plus a proximal regularization term to get an approximation to the original function $F$. Thus, our ALM algorithm can also be viewed as a proximal point algorithm.

A drawback of Algorithm \ref{alg:ALM} is that it requires both $f$ and $g$ to be continuously differentiable. In many applications, however, one of these functions is nonsmooth,  as in the examples given in Section \ref{sec:intro}. Although Algorithm \ref{alg:Aug-Lag-symmetric} can be applied when $f(x)$ and $g(x)$ are nonsmooth, we are unable to provide a comparable complexity bound in this case. However, when only one of the functions of $f$ and $g$ is nonsmooth (say $g$ is nonsmooth), the following  variant of Algorithm \ref{alg:ALM} applies, and for this algorithm, we have a complexity result.

\begin{algorithm2e}\caption{Alternating Linearization Method with Skipping Steps (ALM-S)}
    \label{alg:ALM-S}
\linesnumberedhidden \dontprintsemicolon Choose $\mu$, $\lambda^0$ and
$x^0=y^0$.\; \For{$k=0,1,\cdots$}{
$ x^{k+1}:= \arg\min_{x} \LCal_\mu(x,y^k;\lambda^k)$ \;
{\bf If} $F(x^{k+1})> \LCal_\mu(x^{k+1},y^k;\lambda^k)$, {\bf then} $x^{k+1}  := y^k$ \;
$ y^{k+1}:= \arg\min_{y} Q_f(y,x^{k+1}) $ \;
$\lambda^{k+1} := \nabla f(x^{k+1})-(x^{k+1}-y^{k+1})/\mu$ \;
}
\end{algorithm2e}

We call Algorithm \ref{alg:ALM-S}, ALM with skipping steps (ALM-S) because in line 4 of Algorithm \ref{alg:ALM-S}, if \bea\label{eq:check-ALM-S}F(x^{k+1})> \LCal_\mu(x^{k+1},y^k;\lambda^k)\eea holds, we let $x^{k+1} := y^k$, i.e., we skip the computation of $x^{k+1}$ in line 3. An alternative version of Algorithm \ref{alg:ALM-S} that has smaller average work per iteration is the following Algorithm \ref{alg:ALM-S-equiv}.

\begin{algorithm2e}\caption{Alternating Linearization Method with Skipping Steps (equivalent version)}
    \label{alg:ALM-S-equiv}
\linesnumberedhidden \dontprintsemicolon Choose $\mu$, $\lambda^0$ and
$x^0=y^0$.\; \For{$k=0,1,\cdots$}{
$ x^{k+1}:= \arg\min_{x} \LCal_\mu(x,y^k;\lambda^k)$ \;
\eIf{$F(x^{k+1})> \LCal_\mu(x^{k+1},y^k;\lambda^k)$}{
   $x^{k+1}  := y^k$ \;
   $ y^{k+1}:= \arg\min_{y} Q_f(y,x^{k+1}) $ \;
   $\lambda^{k+1} := \nabla f(x^{k+1})-(x^{k+1}-y^{k+1})/\mu$ \;
   }{
   $\lambda^{k+\half} := \lambda^k - (x^{k+1}-y^k)/\mu$ \;
   $ y^{k+1}:= \arg\min_y \LCal_\mu(x^{k+1},y;\lambda^{k+\half})$ \;
   $\lambda^{k+1} := \lambda^{k+\half} - (x^{k+1} - y^{k+1})/\mu$ \;
   }
}
\end{algorithm2e}

Note that in Algorithm \ref{alg:ALM-S-equiv}, when \eqref{eq:check-ALM-S} does not hold, we switch to the SADAL algorithm, which updates $\lambda$ instead of computing $\nabla f(x^{k+1})$. Algorithm \ref{alg:ALM-S-equiv} is usually faster than Algorithm \ref{alg:ALM-S} when $\nabla f(x^{k+1})$ is costly to compute in addition to performing Step 3. Note also that when \eqref{eq:check-ALM-S} holds, then the steps of the algorithm reduce to those of ISTA.

The following theorem gives conditions under which Algorithms \ref{alg:Aug-Lag-symmetric}, \ref{alg:ALM}, \ref{alg:ALM-S} and \ref{alg:ALM-S-equiv} are equivalent.

\begin{theorem}\label{the:equiv-algs}
(i) If both $f$ and $g$ are differentiable, and $\lambda_0$ is set to $-\nabla g(y^0)$, then Algorithms \ref{alg:Aug-Lag-symmetric} and \ref{alg:ALM} are equivalent.
(ii) If in addition $g$ is Lipschitz continuous with Lipschitz constant $L(g)$, and $\mu\leq 1/L(g)$, then Algorithms \ref{alg:ALM} and \ref{alg:ALM-S} are equivalent.
(iii) If $f$ is differentiable, then Algorithms \ref{alg:ALM-S} and \ref{alg:ALM-S-equiv} are equivalent.
\end{theorem}
\begin{proof}
When both $f$ and $g$ are differentiable and $\lambda_0=-\nabla g(y^0)$, \eqref{gradg,gradf,lambda} holds for all $k\geq 0$, and it follows that $\LCal_\mu(x,y^k;\lambda^k)\equiv Q_g(x,y^k)$ and $\LCal_\mu(x^{k+1},y;\lambda^{k+\half})\equiv Q_f(y,x^{k+1}).$ This proves part (i). If $\nabla g(x)$ is Lipschitz continuous and $\mu\leq 1/L(g)$, \beaa g(x^{k+1}) \leq g(y^k) + \langle \nabla g(y^k), x^{k+1}-y^k \rangle + \oneotwomu\left\|x^{k+1}-y^k\right\|_2^2, \eeaa holds (see e.g., \citep{Bertsekas-book-99}).
This implies that \eqref{eq:check-ALM-S} does not hold and hence, $x^{k+1}:= \arg\min_{x} \LCal_\mu(x,y^k;\lambda^k)$, and the equivalence of Algorithms \ref{alg:ALM} and \ref{alg:ALM-S} follows. This proves part (ii). The optimality of $x^{k+1}$ in line 3 of Algorithm \ref{alg:ALM-S-equiv} implies that $\lambda^{k+\half}=\nabla f(x^{k+1})$ when \eqref{eq:check-ALM-S} does not hold and hence that $\LCal_\mu(x^{k+1},y;\lambda^{k+\half})\equiv Q_f(y,x^{k+1})$. This proves part (iii).
\end{proof}

We show in the following that the iteration complexity
of Algorithm \ref{alg:ALM-S} is $O(1/\epsilon)$ for obtaining an $\epsilon$-optimal solution for \eqref{prob:P}.
First, we need the following generalization of Lemma 2.3 in \citep{Beck-Teboulle-2009}.

\begin{lemma}\label{lem:BT2.3}
Let $\psi: \br^n\rightarrow \br$ and $\phi:\br^n\rightarrow \br$ be convex functions and define
\beaa Q_\psi(u,v):= \phi(u) + \psi(v) + \langle \gamma_\psi(v), u-v \rangle + \frac{1}{2\mu}\|u-v\|_2^2,  \eeaa and \bea\label{def:p_g(v)} p_\psi(v):=\arg\min_u Q_\psi(u,v),\eea
where $\gamma_\psi(v)$ is any subgradient in the subdifferential $\partial\psi(v)$ of $\psi(v)$ at the point $v$.
Let $\Phi(\cdot) = \phi(\cdot)+\psi(\cdot)$. For any  $v$, if
\bea\label{lem:BT2.3-assump-g} \Phi(p_\psi(v)) \leq Q_\psi(p_\psi(v),v), \eea
then for any $u$,
\bea\label{lem:BT2.3-conclude-g} 2\mu(\Phi(u) - \Phi(p_\psi(v))) \geq \|p_\psi(v)-u\|^2 - \|v-u\|^2. \eea
\end{lemma}

{\em Proof.}
From \eqref{lem:BT2.3-assump-g}, we have \bea\label{lem:BT2.3-proof-eq-1}\ba{lll}\Phi(u) - \Phi(p_\psi(v)) & \geq & \Phi(u) - Q_\psi(p_\psi(v),v) \\
      & = & \Phi(u) - \left(\phi(p_\psi(v)) + \psi(v) + \langle \gamma_\psi(v),p_\psi(v)-v\rangle + \frac{1}{2\mu}\|p_\psi(v)-v\|_2^2\right).\ea\eea
Since $\phi$ and $\psi$ are convex we have
\bea\label{lem:BT2.3-proof-eq-a} \phi(u) \geq \phi(p_\psi(v)) + \langle u-p_\psi(v), \gamma_\phi(p_\psi(v)) \rangle, \eea and
\bea\label{lem:BT2.3-proof-eq-b} \psi(u) \geq \psi(v) + \langle u- v, \gamma_\psi(v)\rangle, \eea where $\gamma_\phi(\cdot)$ is a subgradient of $\phi(\cdot)$ and  $\gamma_\phi(p_\psi(v))$ satisfies the first-order optimality conditions for \eqref{def:p_g(v)}, i.e.,
\bea\label{def:p_g(v)-1st-opt-cond}\gamma_\phi(p_\psi(v)) + \gamma_\psi(v) + \frac{1}{\mu}(p_\psi(v)-v)=0.\eea Summing \eqref{lem:BT2.3-proof-eq-a} and \eqref{lem:BT2.3-proof-eq-b} yields
\bea\label{lem:BT2.3-proof-eq-2}\Phi(u) \geq \phi(p_\psi(v)) + \langle u-p_\psi(v), \gamma_\phi(p_\psi(v)) \rangle + \psi(v) + \langle u- v, \gamma_\psi(v)\rangle. \eea
Therefore, from \eqref{lem:BT2.3-proof-eq-1}, \eqref{def:p_g(v)-1st-opt-cond} and \eqref{lem:BT2.3-proof-eq-2} it follows that
\bea\label{lem:BT2.3-proof-eq-3}\begin{split} \Phi(u) - \Phi(p_\psi(v)) & \geq \langle \gamma_\psi(v)+\gamma_\phi(p_\psi(v)), u-p_\psi(v) \rangle - \frac{1}{2\mu}\|p_\psi(v)-v\|_2^2 \\
                                & = \langle -\frac{1}{\mu}(p_\psi(v)-v), u-p_\psi(v)\rangle - \frac{1}{2\mu}\|p_\psi(v)-v\|_2^2 \\
                                & = \frac{1}{2\mu}\left(\|p_\psi(v)-u\|^2-\|v-u\|^2\right).\qquad \endproof\end{split} \eea

\begin{theorem}\label{the:ALM-S}
Assume $\nabla f(\cdot)$  is Lipschitz continuous with Lipschitz constant $L(f)$. For $\mu \leq 1/L(f)$, the iterates $y^k$ in
Algorithm \ref{alg:ALM-S} satisfy
\begin{align}\label{the:ALM-nonsmooth-conclude} F(y^k)-F(x^*)\leq\frac{\|x^0-x^*\|^2}{2\mu (k+k_n)}, \quad\forall k, \end{align} where $x^*$ is an optimal solution of \eqref{prob:P} and $k_n$ is the number of
iterations until the $k$-th for which $F(x^{k+1})\leq \LCal_\mu(x^{k+1}, y^k; \lambda^k)$, i.e., the number of iterations when no skipping step occurs. Thus, the sequence $\{F(y^k)\}$ produced by Algorithm \ref{alg:ALM-S} converges to $F(x^*)$. Moreover, if $1/(\beta L(f)) \leq \mu \leq 1/  L(f)$ where $\beta\geq 1$, the number of iterations needed to obtain an $\epsilon$-optimal solution is at most $\lceil C/\epsilon\rceil$, where $C=\beta L(f)\|x^0-x^*\|^2/2$.
\end{theorem}

\begin{proof}
Let $I$ be the set of all iteration indices until $k-1$-st for which no skipping
occurs and let $I_c$ be its complement. Let $I=\{n_i\}, \ i=0, \ldots, k_n-1$.
It follows that for all $n\in I_c$, $x^{n+1}=y^n$.

For $n\in I$ we can apply Lemma \ref{lem:BT2.3} to obtain
the following inequalities.
In \eqref{lem:BT2.3-conclude-g}, by letting $\psi=f$, $\phi=g$, $u=x^*$ and $v=x^{n+1}$, we
get $p_\psi(v)=y^{n+1}$, $\Phi=F$ and \bea\label{proof-ALM-1}2\mu(F(x^*)-F(y^{n+1}))\geq
\|y^{n+1}-x^*\|^2-\|x^{n+1}-x^*\|^2.\eea Similarly, by letting $\psi=g$, $\phi=f$,
$u=x^*$ and $v=y^{n}$ in \eqref{lem:BT2.3-conclude-g} we get $p_g(v)=x^{n+1}$, $\Phi=F$ and
\bea\label{proof-ALM-2}2\mu(F(x^*)-F(x^{n+1}))\geq
\|x^{n+1}-x^*\|^2-\|y^{n}-x^*\|^2.\eea Taking the summation of
\eqref{proof-ALM-1} and \eqref{proof-ALM-2} we get
\bea\label{proof-ALM-3}2\mu(2F(x^*)-F(x^{n+1})-F(y^{n+1}))\geq
\|y^{n+1}-x^*\|^2-\|y^n-x^*\|^2.\eea

For $n\in I_c$, \eqref{proof-ALM-1} holds. Then since $x^{n+1}=y^n$ we get
\bea\label{proof-ALM-3.5}2\mu(F(x^*)-F(y^{n+1}))\geq
\|y^{n+1}-x^*\|^2-\|y^n-x^*\|^2.\eea

Summing \eqref{proof-ALM-3} and  \eqref{proof-ALM-3.5}
over $n=0,1,\ldots,k-1$ we get \begin{align}\label{proof-ALM-4}
& 2\mu((2|I|+|I_c|)F(x^*)-\sum_{n\in I}F(x^{n+1}) - \sum_{n=0}^{k-1}F(y^{n+1})) \\
\nonumber\geq & \sum_{n=0}^{k-1}\left(\|y^{n+1}-x^*\|^2-\|y^{n}-x^*\|^2\right) \\
\nonumber = & \|y^k-x^*\|^2-\|y^0-x^*\|^2 \\
\nonumber \geq & - \|x^0-x^*\|^2.\end{align}

For any $n$, since Lemma \ref{lem:BT2.3} holds for any $u$, letting $u=x^{n+1}$ instead of $x^*$ we get from \eqref{proof-ALM-1} that
\bea\label{proof-ALM-5}2\mu(F(x^{n+1})-F(y^{n+1}))\geq
\|y^{n+1}-x^{n+1}\|^2 \geq 0,\eea
 or, equivalently,
\bea\label{proof-ALM-6}2\mu(F(x^{n})-F(y^{n}))\geq
\|y^{n}-x^{n}\|^2 \geq 0.\eea
Thus we get $F(y^n)\leq F(x^n), \forall n.$

Similarly, for $n\in I$ by letting
$u=y^{n}$ instead of $x^*$ we get from \eqref{proof-ALM-2} that
\bea\label{proof-ALM-7}2\mu(F(y^{n})-F(x^{n+1}))\geq\|x^{n+1}-y^{n}\|^2 \geq 0.\eea
On the other hand, for $n\in I_c$, \eqref{proof-ALM-7} holds trivially because
$x^{n+1}=y^n$; thus \eqref{proof-ALM-7} holds for all $n$.

Adding \eqref{proof-ALM-5} and \eqref{proof-ALM-7} and adding \eqref{proof-ALM-6} and \eqref{proof-ALM-7}, respectively, yield
\bea\label{proof-ALM-9}2\mu(F(y^n)-F(y^{n+1})) \geq 0 \mbox{ and } 2\mu(F(x^{n})-F(x^{n+1}))\geq 0, \mbox{ for all } n. \eea
The inequalities \eqref{proof-ALM-9} show that the sequences of function values $F(y^n)$ and $F(x^n)$ are non-increasing. Thus we have,
\bea\label{proof-ALM-12} \sum_{n=0}^{k-1}F(y^{n+1})\geq k F(y^k)\quad \mbox{ and } \quad \sum_{n\in I}F(x^{n+1}) \geq k_n F(x^k).\eea

Combining \eqref{proof-ALM-4} and \eqref{proof-ALM-12} yields
\bea\label{proof-ALM-13}2\mu\left((k+k_n)F(x^*)-k_n F(x^k)-k F(y^k)\right)\geq -\|x^0-x^*\|^2.\eea
Hence, since $F(y^k)\leq F(x^k)$, \beaa 2\mu (k+k_n) \left(F(y^k)-F(x^*)\right) \leq \|x^0-x^*\|^2, \eeaa which gives us the desired result \eqref{the:ALM-nonsmooth-conclude}. \end{proof}

\begin{corollary}\label{cor:ALM-nonsmooth}
Assume $\nabla f$ and $\nabla g$ are both Lipschitz continuous with Lipschitz constants $L(f)$ and $L(g)$, respectively. For $\mu \leq \min\{1/L(f),1/L(g)\}$,
Algorithm \ref{alg:ALM} satisfies
\begin{align}\label{cor:ALM-nonsmooth-conclude} F(y^k)-F(x^*)\leq\frac{\|x^0-x^*\|^2}{4\mu k}, \quad \forall k, \end{align} where $x^*$ is an optimal solution of \eqref{prob:P}. Thus sequence $\{F(y^k)\}$ produced by Algorithm \ref{alg:ALM} converges to $F(x^*)$. Moreover, if $ 1/(\beta\max\{L(f),L(g)\})\leq \mu
\leq 1/\max\{L(f),L(g)\}$ where $\beta\geq 1$, the number of iterations needed to get an
 $\epsilon$-optimal solution is at most $\lceil C/\epsilon\rceil$, where $C=\beta \max\{L(f), L(g)\}\|x^0-x^*\|^2/4$.
\end{corollary}
\begin{proof}
The conclusion follows from Theorems \ref{the:equiv-algs} and \ref{the:ALM-S} and $k_n=k$. \end{proof}

\begin{remark}
The complexity bound in Corollary \ref{cor:ALM-nonsmooth} is smaller than the analogous bound for ISTA in \citep{Beck-Teboulle-2009} by a factor of two.
It is easy to see that the bound in Theorem \ref{the:ALM-S} is also an improvement over the bound in \citep{Beck-Teboulle-2009} as long as
$F(x^{k+1})\leq \LCal_\mu(x^{k+1}, y^k; \lambda^k)$ holds for at least
one value of $k$. It is reasonable then to ask if the per-iteration cost of Algorithms \ref{alg:ALM} and \ref{alg:ALM-S} are comparable to that of ISTA.
It is indeed the case when the assumption holds that minimizing $ \LCal_\mu(x, y^k; \lambda^k)$ has comparable cost (and often involves the same computations) as computing the gradient $\nabla f(y^k)$.
\end{remark}
\begin{remark}
More general problems of the form \bea \label{prob:general} \ba{ll} \min & f(x) + g(y) \\ \st & Ax + y =b \ea \eea are easily handled by our approach, since one can express \eqref{prob:general} as \[\min \quad f(x) + g(b-Ax).\]
\end{remark}

\begin{remark}\label{rem:convex-set}
If a convex constraint $x\in\mathcal{C}$, where $\mathcal{C}$ is a convex set is added to problem \eqref{prob:P}, and we impose this constraint in the two  subproblems in Algorithms \ref{alg:ALM} and \ref{alg:ALM-S}, i.e., we impose $x\in\mathcal{C}$ in the subproblems with respect to $x$ and $y\in\mathcal{C}$ in the subproblems with respect to $y$, the complexity results in Theorem \ref{the:ALM-S} and Corollary \ref{cor:ALM-nonsmooth} continue to hold. The only changes in the proof are in Lemma \ref{lem:BT2.3}. If there is a constraint $x\in\mathcal{C}$, then \eqref{lem:BT2.3-conclude-g} holds for any $u\in\CCal$ and $v\in\CCal$. Also in the proof of Lemma \ref{lem:BT2.3}, the first equality in \eqref{lem:BT2.3-proof-eq-3} becomes a ``$\geq$'' inequality due to the fact that the optimality conditions \eqref{def:p_g(v)-1st-opt-cond} become \beaa \langle \gamma_\phi(p_\psi(v)) + \gamma_\psi(v) + \frac{1}{\mu}(p_\psi(v)-v), u - p_\psi(v)\rangle \geq 0, \forall u\in\mathcal{C}.\eeaa
\end{remark}
\begin{remark}
Although Algorithms \ref{alg:ALM} and \ref{alg:ALM-S} assume that the Lipschitz constants are
known, and hence that an upper bound for $\mu$ is known, this can be relaxed by using
the backtracking technique in \citep{Beck-Teboulle-2009} to estimate
$\mu$ at each iteration.
\end{remark}

\section{Fast Alternating Linearization Methods}\label{sec:FALM} In this section, we propose a fast alternating
linearization method (FALM) which computes an $\epsilon$-optimal solution to problem \eqref{prob:P} in $O(\sqrt{L/\epsilon})$ iterations, while keeping
the work at each iteration almost the same as that required by ALM.

FALM is an accelerated version of ALM for solving \eqref{prob:P}, or
equivalently \eqref{prob:P-x-y}, when $f(x)$ and $g(x)$ are both differentiable, and is given below as
Algorithm \ref{alg:FALM}. Clearly, FALM is also a Gauss-Seidel type algorithm. In fact, it is a successive over-relaxation type algorithm since $(t_k-1)/t_{k+1}>0, \forall k\geq 2.$

\begin{algorithm2e}\caption{Fast Alternating Linearization Method (FALM)}
    \label{alg:FALM}
\linesnumberedhidden \dontprintsemicolon Choose $\mu$ and $x^0=y^0=z^1$, set $t_1=1$.\; \For{$k=1,2,\cdots$}{
$x^{k} :=  \arg\min_x Q_g(x,z^k)$ \;
$y^{k} :=  \arg\min_y Q_f(y,x^{k})$ \;
$t_{k+1} := (1+\sqrt{1+4t_k^2})/2$ \;
$z^{k+1} := y^k + \frac{t_k-1}{t_{k+1}}(y^k-y^{k-1})$ \;
}
\end{algorithm2e}

Algorithm \ref{alg:FALM} requires both $f$ and $g$ to be continuously differentiable. To develop an algorithm that can be applied to problems where one of the functions is non-differentiable, we use a skipping technique as in Algorithm \ref{alg:ALM-S}. FALM with skipping steps (FALM-S), which does not require $g(x)$ to be smooth, is given below as Algorithm \ref{alg:FALM-S}.
\begin{algorithm2e}\caption{FALM with Skipping Steps (FALM-S)} \label{alg:FALM-S}
\linesnumberedhidden \dontprintsemicolon
Choose $x^0=y^0=z^1$ and $\lambda^1\in-\partial g(z^1)$, set $t_1=1$.\; \For{$k=1,2,\cdots$}{
$x^{k} := \arg\min_x \LCal_\mu(x,z^k;\lambda^k)$ \;
\If {$F(x^k) > \LCal_\mu(x^k,z^k;\lambda^k)$} {\eIf {x-step was not skipped at iteration k-1}{$t_k:=\left(1+\sqrt{1+8t_{k-1}^2}\right)/2$}{$t_k:=\left(1+\sqrt{1+4t_{k-1}^2}\right)/2$}\;$x^k:=z^k:=y^{k-1}+\frac{t_{k-1}-1}{t_k}\left(y^{k-1}-y^{k-2}\right)$}  \;
$y^{k} := \arg\min_y Q_f(x^{k},y)$ \;
\eIf {$x^k=z^k$}{$t_{k+1}:= \left(1+\sqrt{1+2t_k^2}\right)/2$}{$t_{k+1}:= \left(1+\sqrt{1+4t_k^2}\right)/2$} \;
$z^{k+1} := y^k + \frac{t_k-1}{t_{k+1}}\left(y^k-y^{k-1}\right)$ \;
Choose $\lambda^{k+1}\in -\partial g(z^{k+1})$ \;
}
\end{algorithm2e}

The following theorem gives conditions under which Algorithms \ref{alg:FALM} and \ref{alg:FALM-S} are equivalent.
\begin{theorem}\label{the:equiv-algs-FALM}
If both $f(x)$ and $g(x)$ are differentiable and $\nabla g(x)$ is Lipschitz continuous with Lipschitz constant $L(g)$, and $\mu\leq 1/L(g)$, then Algorithms \ref{alg:FALM} and \ref{alg:FALM-S} are equivalent.
\end{theorem}
\begin{proof}
As in Theorem \ref{the:equiv-algs}, if $f$ and $g$ are differentiable, $\LCal_\mu(x,z^k;\lambda^k)\equiv Q_g(x,z^k)$. The conclusion then follows from the fact that $F(x^k) \leq \LCal_\mu(x^k,z^k;\lambda^k)$
always holds when $\mu\leq 1/L(g)$ and thus there are no skipping steps.
\end{proof}

To prove that Algorithm \ref{alg:FALM-S} requires $O(\sqrt{L(f)/\epsilon})$ iterations to obtain an $\epsilon$-optimal solution, we need the following lemmas. We call $k$-th iteration a {\it skipping step} if $x^k=z^k$, and a {\it regular step} if $x^k\neq z^k$.
\begin{lemma}\label{lem:FALM-S}
The sequence $\{x^k,y^k\}$ generated by Algorithm \ref{alg:FALM-S} satisfies \bea\label{lem:FALM-S-eq} 2\mu(t_k^2v_k - t_{k+1}^2v_{k+1}) \geq \|u^{k+1}\|^2 - \|u^k\|^2, \eea where $u^k:=t_ky^k-(t_k-1)y^{k-1}-x^*$ and $v_k:=2F(y^k)-2F(x^*)$ if iteration $k$ is a regular step and $v_k:=F(y^k)-F(x^*)$ if iteration $k$ is a skipping step.
\end{lemma}

\begin{proof}
There are four cases to consider: (i) both the $k$-th and the $(k+1)$-st iterations are regular steps; (ii) the $k$-th iteration is a regular step and the $(k+1)$-st iteration is a skipping step; (iii) both the $k$-th and the $(k+1)$-st iterations are skipping steps; (iv) the $k$-th iteration is a skipping step and the $(k+1)$-st iteration is a regular step. We will prove that the following inequality holds for all the four cases:
\bea\label{proof-FALM-tk^2vk-tk+1^2vk+1}\begin{split} & 2\mu(t_k^2v_k-t_{k+1}^2v_{k+1}) \\ \geq &
t_{k+1}(t_{k+1}-1)\left(\|y^{k+1}-y^k\|^2-\|z^{k+1}-y^k\|^2\right)+t_{k+1}\left(\|y^{k+1}-x^*\|^2-\|z^{k+1}-x^*\|^2\right).\end{split}\eea
The proof of \eqref{lem:FALM-S-eq} and hence, the lemma, then follows from the fact that the right hand side of inequality \eqref{proof-FALM-tk^2vk-tk+1^2vk+1} equals
\beaa \|t_{k+1}y^{k+1}-(t_{k+1}-1)y^k-x^*\|^2-\|t_{k+1}z^{k+1}-(t_{k+1}-1)y^k-x^*\|^2 = \|u^{k+1}\|^2-\|u^k\|^2, \eeaa where we have used the fact that $t_{k+1}z^{k+1}:=t_{k+1}y^k+t_k(y^k-y^{k-1})-(y^k-y^{k-1}).$

Case (i):
Let us first consider the case when both the $k$-th and $(k+1)$-st iterations are regular steps.
In \eqref{lem:BT2.3-conclude-g}, by letting $\psi=f$, $\phi=g$, $u=y^k$ and
$v=x^{k+1}$, we get $p_\psi(v) =y^{k+1}$, $\Phi=F$ and
\bea\label{proof-FALM-S-1} 2\mu(F(y^k) - F(y^{k+1}))\geq \|y^{k+1}-y^k\|^2-\|x^{k+1}-y^k\|^2.\eea
In \eqref{lem:BT2.3-conclude-g}, by letting $\psi=g$, $\phi=f$, $u=y^k$, $v=z^{k+1}$, we get $p_\psi(v)=x^{k+1}$, $\Phi=F$ and \bea\label{proof-FALM-S-2} 2\mu(F(y^k) - F(x^{k+1}))\geq \|x^{k+1}-y^k\|^2-\|z^{k+1}-y^k\|^2.\eea
Summing \eqref{proof-FALM-S-1} and \eqref{proof-FALM-S-2}, and
using the fact that $F(y^{k+1})\leq F(x^{k+1})$, we obtain, \bea\label{proof-FALM-vk-vk+1}
2\mu(v_k-v_{k+1}) = 2\mu(2F(y^k)-2F(y^{k+1})) \geq
\|y^{k+1}- y^k\|^2-\|z^{k+1}- y^k\|^2\eea Again, in
\eqref{lem:BT2.3-conclude-g}, by letting $\psi=g$, $\phi=f$, $u=x^*$, $v=z^{k+1}$, we get $p_\psi(v)=x^{k+1}$, $\Phi=F$ and
\bea\label{proof-FALM-S-3}2\mu(F(x^*)-F(x^{k+1}))\geq\|x^{k+1}-x^*\|^2-\|z^{k+1}-x^*\|^2.\eea
In \eqref{lem:BT2.3-conclude-g}, by letting $\psi=f$, $\phi=g$, $u=x^*$, $v=x^{k+1}$, we get $p_\psi(v)=y^{k+1}$, $\Phi=F$ and
\bea\label{proof-FALM-S-4}2\mu(F(x^*)-F(y^{k+1})\geq\|y^{k+1}-x^*\|^2-\|x^{k+1}-x^*\|^2.\eea
Summing \eqref{proof-FALM-S-3} and
\eqref{proof-FALM-S-4}, and again using the fact that $F(y^{k+1})\leq F(x^{k+1})$, we obtain, \bea\label{proof-FALM--vk+1}-2\mu
v_{k+1} = 2\mu (2F(x^*)-2F(y^{k+1})) \geq \|y^{k+1}-x^*\|^2-\|z^{k+1}-x^*\|^2.\eea If we multiply \eqref{proof-FALM-vk-vk+1} by $t_k^2$,
and \eqref{proof-FALM--vk+1} by $t_{k+1}$, and take the sum of the
resulting two inequalities, we get \eqref{proof-FALM-tk^2vk-tk+1^2vk+1} by using
the fact that $t_k^2=t_{k+1}(t_{k+1}-1)$.

Case (ii): By letting $\psi=f$, $\phi=g$, $u=y^k$ and $v=z^{k+1}$ in \eqref{lem:BT2.3-conclude-g}, we get $p_\psi(v)=y^{k+1}$, $\Phi=F$ and
\bea\label{proof-FALM-N-1-skip} 2\mu(F(y^k) - F(y^{k+1}))\geq \|y^{k+1}-y^k\|^2-\|z^{k+1}- y^k\|^2.\eea
Since the steps taken in the $k$-th and $(k+1)$-st iterations are regular and skipping, respectively, we have
\bea\label{proof-FALM-vk-vk+1-skip}
2\mu\left(\frac{v_k}{2}-v_{k+1}\right) = 2\mu(F(y^k)-F(y^{k+1}))\geq
\|y^{k+1}-y^k\|^2-\|z^{k+1}-y^k\|^2\eea
Also by letting $\psi=f$, $\phi=g$, $u=x^*$ and $v=z^{k+1}$ in \eqref{lem:BT2.3-conclude-g}, we get $p_\psi(v)=y^{k+1}$, $\Phi=F$ and
\bea\label{proof-FALM-N-4-skip}
-2\mu v_{k+1} = 2\mu(F(x^*)-F(y^{k+1}))\geq\|y^{k+1}-x^*\|^2-\|z^{k+1}-x^*\|^2.\eea
Then multiplying \eqref{proof-FALM-vk-vk+1-skip} by $2t_k^2$,
\eqref{proof-FALM-N-4-skip} by $t_{k+1}$, summing the
resulting two inequalities and using
the fact that in this case $2t_k^2=t_{k+1}(t_{k+1}-1)$, we obtain \eqref{proof-FALM-tk^2vk-tk+1^2vk+1}.

Case (iii): This case reduces to two consecutive FISTA steps and the proof above applies with $t_k^2=t_{k+1}(t_{k+1}-1)$ and inequality \eqref{proof-FALM-vk-vk+1-skip} replaced by
\bea\label{proof-FALM-vk-vk+1-skipskip}
2\mu(v_k-v_{k+1})=2\mu(F(y^k)-F(y^{k+1}))\geq
\|y^{k+1}-y^k\|^2-\|z^{k+1}-y^k\|^2\eea
which gets multiplied by $t_k^2$.

Case (iv): In this case, \eqref{proof-FALM-vk-vk+1} in the proof of case (i) is replaced by
\bea\label{proof-FALM-vk-vk+1-skipnoskip}
2\mu(2v_k-v_{k+1})=2\mu(2F(y^k)-2F(y^{k+1}))\geq
\|y^{k+1}-y^k\|^2-\|z^{k+1}-y^k\|^2\eea
which when multiplied by $t_k^2/2$ and combined with \eqref{proof-FALM--vk+1} multiplied by $t_{k+1}$, and the fact that in this case
$t_k^2/2=t_{k+1}(t_{k+1}-1)$, yields \eqref{proof-FALM-tk^2vk-tk+1^2vk+1}.
\end{proof}

The following lemma gives lower bounds for the sequence of scalars $\{t_k\}$ generated by Algorithm \ref{alg:FALM-S}.
\begin{lemma}\label{lem:tk-bound-FALM-S}
For all $k \geq 1$ the sequence $\{t_k\}$ generated by Algorithm \ref{alg:FALM-S} satisfies:

\noindent
if the first step is a skipping step,
\beaa
t_k \geq
\left\{\begin{array}{lll}
\frac{1}{2}(k + 1 + \alpha r(k)) & \mbox{if} & k \quad \mbox{is a skipping step}, \\
\frac{1}{2\sqrt{2}}(k + 1 + \alpha r(k)) & \mbox{if} & k \quad \mbox{is a regular step},
\end{array}\right.
\eeaa
if the first step is a regular step,
\beaa
t_k \geq
\left\{\begin{array}{lll}
\frac{1}{\sqrt{2}}(k + 1 + \hat{\alpha} s(k)) & \mbox{if} & k \quad \mbox{is a skipping step},\\
\frac{1}{2}(k + 1 + \hat{\alpha} s(k)) & \mbox{if} & k \quad \mbox{is a regular step},
\end{array}\right.
\eeaa
where $r(k)$ and $s(k)$ are the number of steps among the first $k$ steps that are regular and skipping  steps, respectively, and $\alpha \equiv \sqrt{2} - 1$ and $\hat{\alpha} \equiv \frac{1}{\sqrt{2}} - 1$.
\end{lemma}
\begin{proof}
Consider the case where the first iteration of Algorithm \ref{alg:FALM-S} is a skipping step. Clearly, the sequence of iterations follows a pattern of alternating blocks of one or more skipping steps and one or more regular steps. Let the index of the first iteration in the $i$-th block be denoted by $n_i$.
Since it is assumed that the first iteration is a skipping step, iterations $n_1,n_3,n_5\ldots$ are skipping steps ($n_1 = 1$) and $n_2,n_4,n_6\ldots$ are regular steps. Note that the statement of the lemma in this case corresponds to
\bea\label{lem:tk-bd}
t_k \geq
\left\{\begin{array}{ll}
\frac{1}{2}(k + 1 + \alpha r(k)), & \mbox{for $n_j \leq k \leq n_{j+1}-1$, if $j$ is odd}, \\
\frac{1}{2\sqrt{2}}(k + 1 + \alpha r(k)), & \mbox{for $n_j \leq k \leq n_{j+1}-1$, if $j$ is even},
\end{array}\right.
\eea
which we will prove by induction on $j$.

We first note that it follows from the updating rules and formulas for $t_k$ that

\bea\label{lem:tk}
t_k \geq
\left\{\begin{array}{ll}
\frac{1}{2} + \sqrt{2} t_{k-1}, & \mbox{if } k  \mbox{ is a skipping step and } k-1 \mbox{ is a regular step}, \\
\frac{1}{2} + \frac{1}{\sqrt{2}}t_{k-1}, & \mbox{if }  k  \mbox{ is a regular step and } k-1 \mbox{ is a skipping step}, \\
\frac{1}{2} + t_{k-1}, & \mbox{otherwise}.
\end{array}\right.
\eea
Consider $j = 1$. Clearly (\ref{lem:tk-bd}) holds for all iterations $n_1 = 1 \leq k \leq n_2 - 1$,
since $t_k\geq \frac{k+1}{2}$ holds trivially for $t_1=1$, and for $1 < k \leq n_2-1$, $t_{k}\geq \frac{k-1}{2}+t_1 = \frac{k+1}{2}$.

Now assume that (\ref{lem:tk-bd}) holds for all $ j < \bar j$.
If $\bar j$ is even, iterations $k = n_{\bar j}$ and $k-1$ are, respectively,  regular and skipping iterations. Hence, from (\ref{lem:tk}), we have that
\[t_{k}\geq\half+\frac{1}{\sqrt{2}}t_{k -1}
\geq \half+\frac{1}{2\sqrt{2}}(k + \alpha r(k-1))=
\frac{1}{2\sqrt{2}}(k + 1 + \alpha r(k)).\]
Since the remaining
$p \equiv n_{{\bar j}+1}-1 - n_{\bar j}$ iterations before iteration
$n_{{\bar j}+1}$ are all regular iterations ($p$ may be zero), we have from (\ref{lem:tk})  that, for $n_{\bar j} < k \leq n_{\bar j+1}-1$,
\beaa \ba{lll} t_{k} & \geq & \frac{k-n_{\bar j}}{2}+ t_{n_{\bar j}} \geq \frac{k-n_{\bar j}}{2}+ \frac{1}{2\sqrt{2}} (n_{\bar j} +1 +\alpha r(n_{\bar j})) \\ & = & \frac{k-n_{\bar j}}{2}+ \frac{1}{2\sqrt{2}} (n_{\bar j} +1 +\alpha (r(k)-k+n_{\bar j}))=
\frac{1}{2\sqrt{2}}(k + 1 + \alpha r(k)).\ea\eeaa

If $\bar j$ is odd, iteration $k = n_{\bar j}$ is a skipping iteration. Hence, from (\ref{lem:tk}), we have that
$t_{k}\geq\half+\sqrt{2}t_{k -1} \geq\half+\sqrt{2}\frac{1}{2\sqrt{2}}(k+\alpha r(k-1)) = \frac{1}{2}(k + 1 + \alpha r(k))$. Since the remaining $p \equiv n_{{\bar j}+1}-1 - n_{\bar j}$ iterations before iteration $n_{{\bar j}+1}$ are all skipping iterations (again $p$ may be zero), we have from (\ref{lem:tk}) that, for $n_{\bar j} < k \leq n_{\bar j+1}-1$,
$t_{k}\geq\frac{k-n_{\bar j}}{2}+ t_{n_{\bar j}} \geq\frac{k-n_{\bar j}}{2}+ \frac{1}{2} (n_{\bar j} +1 +\alpha r(n_{\bar j}))=
\frac{1}{2}(k + 1 + \alpha r(k))$.
This concludes the induction.

Since the proof for the case that the first step is a regular step is totally analogous to the above proof, we leave this to the reader.
\end{proof}

Now we are ready to give the complexity of Algorithm \ref{alg:FALM-S}.
\begin{theorem}\label{the:FALM-S}
Let $\alpha=\sqrt{2}-1$ and $r(k)$ be the number of steps among the first $k$ steps that are regular steps. Assuming $\nabla f(\cdot)$ is Lipschitz continuous with Lipschitz constant $L(f)$, if $\mu\leq 1/L(f)$, the sequence $\{y^k\}$ generated by Algorithm \ref{alg:FALM-S} satisfies:

\be \label{proof-the-FALM-S-bound-v3} F(y^k)-F(x^*)\leq\frac{2\|x^0-x^*\|^2}{\mu(k+1+\alpha \hat{r}(k))^2},\ee
where $\hat{r}(k)=r(k)$ if the first step is a skipping step, and $\hat{r}(k)=r(k)+1$ if the first step is a regular step.

Hence, the sequence $\{F(y^k)\}$ produced by Algorithm \ref{alg:ALM-S} converges to $F(x^*)$. Moreover, if $1/(\beta L(f)) \leq \mu \leq 1/  L(f)$ where $\beta\geq 1$, the number of iterations required by Algorithm \ref{alg:FALM-S} to get an $\epsilon$-optimal solution to \eqref{prob:P} is at most $\lfloor
\sqrt{C/\epsilon}\rfloor$, where $C=2\beta L(f)\|x^0-x^*\|^2$.
\end{theorem}

\begin{proof}
Using the same notation as in Lemmas \ref{lem:FALM-S} and \ref{lem:tk-bound-FALM-S},
\eqref{proof-FALM--vk+1} and \eqref{proof-FALM-N-4-skip} imply that
\beaa -2\mu v_1 \geq \|y^1-x^*\|^2 - \|z^1-x^*\|^2 \eeaa holds whether the first iteration is a skipping step or not. Thus we have
\bea\label{proof-the-FALM-S-bound-v1} 2\mu v_1 + \|y^1-x^*\|^2 \leq \|z^1-x^*\|^2 = \|x^0-x^*\|^2 .\eea From Lemma \ref{lem:FALM-S} we know that the sequence $\{2\mu t_k^2v_k+\|u^k\|^2\}$ is non-increasing. Therefore, we have
\bea \label{proof-the-FALM-S-bound-v2} 2\mu t_k^2 v_k \leq 2\mu t_k^2v_k+\|u^k\|^2 \leq 2\mu t_1^2v_1+\|u^1\|^2 = 2\mu v_1 + \|y^1-x^*\|^2 \leq \|x^0-x^*\|^2, \eea where the equality follows from the facts that $t_1=1$ and $u^1=y^1-x^*$, and the last inequality is from \eqref{proof-the-FALM-S-bound-v1}.

Recall the definition of $v_k$ in Lemma \ref{lem:FALM-S} and the bounds of $t_k$ in Lemma \ref{lem:tk-bound-FALM-S}.  We get (i) and (ii) from \eqref{proof-the-FALM-S-bound-v2}.
  Keeping in mind that $v_k$ has a different expression depending on whether the $k$-th step is a  skipping or a regular step, it follows that the sequence $\{y^k\}$ generated by Algorithm \ref{alg:FALM-S} satisfies:

(i) if the first step is a skipping step, then
\[F(y^k)-F(x^*)\leq\frac{2\|x^0-x^*\|^2}{\mu(k+1+\alpha r(k))^2};\]

(ii) if the first step is a regular step, then
\[F(y^k)-F(x^*)\leq\frac{\|x^0-x^*\|^2}{\mu(k+1+\hat{\alpha} s(k))^2}.\]

It is easy to check that these bounds are equivalent to \eqref{proof-the-FALM-S-bound-v3}, and that the worst case bound on the number of iterations follows from \eqref{proof-the-FALM-S-bound-v3}.
\end{proof}

\begin{corollary}\label{cor:FALM-nonsmooth}
Assume $\nabla f$ and $\nabla g$ are both Lipschitz continuous with Lipschitz constants $L(f)$ and $L(g)$, respectively. For $\mu \leq \min\{1/L(f),1/L(g)\}$,
Algorithm \ref{alg:FALM} satisfies
\begin{align}\label{cor:FALM-nonsmooth-conclude} F(y^k)-F(x^*)\leq\frac{\|x^0-x^*\|^2}{\mu (k+1)^2}, \quad \forall k, \end{align} where $x^*$ is an optimal solution of \eqref{prob:P}. Hence, the sequence $\{F(y^k)\}$ produced by Algorithm \ref{alg:FALM} converges to $F(x^*)$, and if $1/(\beta\max\{L(f),L(g)\})\leq \mu\leq 1/\max\{L(f),L(g)\}$, where $\beta\geq 1$, the number of iterations needed to get an $\epsilon$-optimal solution is at most $\lceil \sqrt{C/\epsilon}-1\rceil$, where $C=\beta \max\{L(f), L(g)\}\|x^0-x^*\|^2$.
\end{corollary}
\begin{proof}
Note that since $\mu\leq\min\{1/L(f),1/L(g)\}$, from Theorem \ref{the:equiv-algs-FALM} we know that Algorithms \ref{alg:FALM} and \ref{alg:FALM-S} are equivalent. That is, every step in Algorithm \ref{alg:FALM-S} is a regular step. Therefore, case (ii) in Theorem \ref{the:FALM-S} holds and $s(k)=0$, which leads to \eqref{cor:FALM-nonsmooth-conclude}. \end{proof}

\begin{remark}
The complexity bound in Corollary \ref{cor:FALM-nonsmooth} is smaller than the analogous bound for FISTA in \citep{Beck-Teboulle-2009} by a factor of $\sqrt{2}$.
It is easy to see that the bound in Theorem \ref{the:FALM-S} is also an improvement over the bound in \citep{Beck-Teboulle-2009} as long as
$F(x^{k+1})\leq \LCal_\mu(x^{k+1}, y^k; \lambda^k)$ holds for at least
one value of $k$. As in the case of Algorithms \ref{alg:ALM} and \ref{alg:ALM-S}, the per-iteration cost of Algorithms \ref{alg:FALM} and \ref{alg:FALM-S} are comparable to that of FISTA.
\end{remark}

\begin{remark}
Line 6 in Algorithm \ref{alg:FALM} and Line 15 in Algorithm \ref{alg:FALM-S} can be changed to:
\[z^{k+1} := w^k + \frac{1}{t_{k+1}}[t_k(y^k-w^{k-1})-(w^k-w^{k-1})],\] where $w^k:=\alpha x^k + (1-\alpha)y^k, \alpha\in(0,1)$, and Theorem \ref{the:FALM-S} and Corollary \ref{cor:FALM-nonsmooth} still hold.
\end{remark}

\begin{remark}
Although Algorithms \ref{alg:FALM} and \ref{alg:FALM-S}, assume that the Lipschitz constants are
known, and hence that an upper bound for $\mu$ is known, this can be relaxed by using
the backtracking technique in \citep{Beck-Teboulle-2009} to estimate
$\mu$ at each iteration.
\end{remark}

\section{Comparison of ALM, FALM, ISTA, FISTA, SADAL and SALSA}\label{sec:comparison}

In this section we compare the performance of our basic and fast ALMs, with and without skipping steps, against ISTA, FISTA, SADAL (Algorithm \ref{alg:Aug-Lag-symmetric}) and an alternating direction augmented Lagrangian method SALSA described in \citep{Afonso-BD-Figueiredo-2009} on a benchmark wavelet-based image deblurring problem from \citep{Figueiredo-Nowak-03}. In this problem, the original image is the well-known Cameraman image of size $256\times 256$ and the observed image is obtained after imposing a uniform blur of size $9\times 9$ (denoted by the operator $R$) and Gaussian noise (generated by the function {\it randn} in MATLAB with a seed of 0 and a standard deviation of $0.56$). Since the coefficient of the wavelet transform of the image is sparse in this problem, one can try to reconstruct the image $u$ from the observed image $b$ by solving the problem:
\bea\label{prob:image-deconvolution}\bar{x}: =\arg\min_x \quad \half\|Ax-b\|_2^2 + \rho\|x\|_1, \eea and setting $u:=W\bar{x}$, where $A:=RW$ and $W$ is the inverse discrete Haar wavelet transform with four levels.
By defining $f(x):=\half\|Ax-b\|_2^2$ and $g(x):=\rho\|x\|_1$, it is clear that \eqref{prob:image-deconvolution} can be expressed in the form of \eqref{prob:P} and can be solved by ALM-S (Algorithm \ref{alg:ALM-S}), FALM-S (Algorithm \ref{alg:FALM-S}), ISTA, FISTA, SALSA and SADAL (Algorithm \ref{alg:Aug-Lag-symmetric}). However, in order to use ALM (Algorithm \ref{alg:ALM}) and FALM (Algorithm \ref{alg:FALM}), we need to smooth $g(x)$ first, since these two algorithms require both $f$ and $g$ to be smooth. Here we apply the smoothing technique introduced by Nesterov \citep{Nesterov-2005} since this technique guarantees that the gradient of the smoothed function is Lipschitz continuous. A smoothed approximation to the $\ell_1$ function $g(x):=\rho\|x\|_1$ with smoothness parameter $\sigma>0$ is \bea\label{def:g-sigma-smooth} g_\sigma(x):=\max\{\langle x,z\rangle - \frac{\sigma}{2}\|z\|_2^2 : \|z\|_\infty\leq\rho\}.\eea

It is easy to show that the optimal solution $z_\sigma(x)$ of \eqref{def:g-sigma-smooth} is  \bea\label{def:Z-sigma(Y)} z_\sigma(x) =\min\{\rho,\max\{x/\sigma,-\rho\}\}.\eea According to Theorem 1 in \citep{Nesterov-2005}, the gradient of $g_\sigma$ is given by $\nabla g_\sigma(x)=z_\sigma(x)$ and is Lipschitz continuous with Lipschitz constant $L(g_{\sigma})=1/\sigma$. After smoothing $g$, we can apply Algorithms \ref{alg:ALM} and \ref{alg:FALM} to solve the smoothed problem: \bea\label{prob:image-deconvolution-sigma}\min_x f(x)+g_\sigma(x).\eea We have the following theorem about the $\epsilon$-optimal solutions of problems \eqref{prob:image-deconvolution} and \eqref{prob:image-deconvolution-sigma}.

\begin{theorem}\label{the:image-deconvolution-epsilon-optimal}
Let $\sigma=\frac{\epsilon}{n\rho^2}$ and $\epsilon>0$. If $x(\sigma)$ is an $\epsilon/2$-optimal solution to \eqref{prob:image-deconvolution-sigma}, then $x(\sigma)$ is an $\epsilon$-optimal solution to \eqref{prob:image-deconvolution}.
\end{theorem}

{\em Proof.} Let $D_g:=\max\{\half\|z\|_2^2 : \|z\|_\infty\leq\rho\}=\half n\rho^2$ and $x^*$ and $x^*(\sigma)$ be optimal solution to problems \eqref{prob:image-deconvolution} and \eqref{prob:image-deconvolution-sigma}, respectively. Note that \bea\label{g-sigma-g-bound-image}g_\sigma(x)\leq g(x)\leq g_\sigma(x)+\sigma D_g,\forall x\in\br^n.\eea Using the inequalities in \eqref{g-sigma-g-bound-image} and the facts that $x(\sigma)$ is an $\epsilon/2$-optimal solution to \eqref{prob:image-deconvolution-sigma} and $\sigma D_g=\frac{\epsilon}{2}$, we have
\begin{align*} f(x(\sigma))+g(x(\sigma))-f(x^*)-g(x^*) & \leq  f(x(\sigma))+g_\sigma(x(\sigma))+\sigma D_g-f(x^*)-g_\sigma(x^*) \\ & \leq f(x(\sigma))+g_\sigma(x(\sigma))+\sigma D_g-f(x^*(\sigma))-g_\sigma(x^*(\sigma)) \\ & \leq \epsilon/2 + \sigma D_g = \epsilon. \qquad \endproof \end{align*}

Thus, to find an $\epsilon$-optimal solution to \eqref{prob:image-deconvolution}, we can apply Algorithms \ref{alg:ALM} and \ref{alg:FALM} to find an $\epsilon/2$-optimal solution to \eqref{prob:image-deconvolution-sigma} with $\sigma=\frac{\epsilon}{n\rho^2}$. The iteration complexity results in Corollaries \ref{cor:ALM-nonsmooth} and \ref{cor:FALM-nonsmooth} hold since the gradient of $g_\sigma$ is Lipschitz continuous. However, the numbers of iterations needed by ALM and FALM to obtain an $\epsilon$-optimal solution to \eqref{prob:image-deconvolution} become $O(1/\epsilon^2)$ and $O(1/\epsilon)$, respectively, due to the fact that the Lipschitz constant $L(g_\sigma)=1/\sigma= \frac{n\rho^2}{\epsilon}=O(1/\epsilon)$.

When Algorithms \ref{alg:ALM-S} and \ref{alg:FALM-S} are applied to solve \eqref{prob:image-deconvolution}, the subproblems \eqref{prob:f-shrink} and \eqref{prob:g-shrink} are easy to solve. Specifically, \eqref{prob:f-shrink} corresponds to solving a linear system which is particularly easy to do because of the special structures of $R$ and $W$ (see \citep{Afonso-BD-Figueiredo-2009,Afonso-BD-Figueiredo-2009b}); \eqref{prob:g-shrink} corresponds to a vector shrinkage operation. When Algorithms \ref{alg:ALM} and \ref{alg:FALM} are applied to solve \eqref{prob:image-deconvolution-sigma}, \eqref{prob:g-shrink} with $g$ replaced by $g_\sigma$ is also easy to solve; its optimal solution is $$x:=z - \tau \min\{\rho, \max\{-\rho,\frac{z}{\tau+\sigma}\}\}.$$


Since ALM is equivalent to SADAL when both functions are smooth, we implemented ALM as SADAL when we solved \eqref{prob:image-deconvolution-sigma}. We also applied SADAL to the nonsmooth problem \eqref{prob:image-deconvolution}. We also implemented ALM-S as Algorithm \ref{alg:ALM-S-equiv} since the latter was usually faster. In all algorithms, we set the initial points $x^0=y^0=\mathbf{0}$, and in FALM and FALM-S we set $z^1=\mathbf{0}$. MATLAB codes for SALSA, FISTA and ISTA (modified from FISTA) were downloaded from http://cascais.lx.it.pt/$\sim$mafonso/salsa.html and their default inputs were used. Moreover, $\lambda^0$ was set to $\mathbf{0}$ in algorithms \ref{alg:Aug-Lag-symmetric}, \ref{alg:ALM-S}, \ref{alg:ALM-S-equiv} and \ref{alg:FALM-S} since $\nabla g_\sigma(x^0) = \mathbf{0}$ and $\mathbf{0}\in -\partial g(x^0)$ when $x^0=\mathbf{0}$. Also, whenever $g(x)$ was smoothed, we set $\sigma = 10^{-6}$. $\mu$ was set to $1$ in all the algorithms since the Lipschitz constant of the gradient of function $\half\|RW(\cdot)-b\|_2^2$ was known to be $1$. We set $\mu$ to 1 even for the smoothed problems. Although this violates the requirement $\mu \leq \frac{1}{L(g_\sigma)}$ in Corollaries \ref{cor:ALM-nonsmooth} and \ref{cor:FALM-nonsmooth}, we see from our numerical results reported below that ALM and FALM still work very well. All of the algorithms tested were terminated after 1000 iterations. The (nonsmoothed) objective function values in \eqref{prob:image-deconvolution} produced by these algorithms at iterations: 10, 50, 100, 200, 500, 800 and 1000 for different choices of $\rho$ are presented in Tables \ref{tab:comp-SALSA-2} and \ref{tab:comp-SALSA-3}. The CPU times (in seconds) and the number of iterations required to reduce the objective function value to below $1.04e+5$ and $8.60e+5$ are reported respectively in the last columns of Tables \ref{tab:comp-SALSA-2} and \ref{tab:comp-SALSA-3}.

All of our codes were written in MATLAB and run in MATLAB 7.3.0 on a Dell Precision 670 workstation with an Intel Xeon(TM) 3.4GHZ CPU and 6GB of RAM.

\begin{table}[ht]{\scriptsize
\begin{center}\caption{Comparison of the algorithms for solving \eqref{prob:image-deconvolution} with $\rho=0.01$}\label{tab:comp-SALSA-2} \scriptsize
\begin{tabular}{|c|c|c|c|c|c|c|c|c|}\hline
solver & \multicolumn{7}{|c|}{obj in k-th iteration} & cpu (iter) \\\hline

       & 10 & 50 & 100 & 200 & 500 & 800 & 1000 & \\\hline

FALM-S & 1.767239e+5  & 1.040919e+5 & 1.004322e+5 & 9.726599e+4 & 9.341282e+4 & 9.182962e+4 & 9.121742e+4 & 24.3 (51) \\\hline 

FALM   & 1.767249e+5  & 1.040955e+5 & 9.899843e+4 & 9.516208e+4 & 9.186355e+4 & 9.073086e+4 & 9.028790e+4 & 23.1 (51) \\\hline 

FISTA  & 1.723109e+5  & 1.061116e+5 & 1.016385e+5 & 9.752858e+4 & 9.372093e+4 & 9.233719e+4 & 9.178455e+4 & 26.0 (69) \\\hline 

ALM-S  & 4.218082e+5  & 1.439742e+5 & 1.212865e+5 & 1.107103e+5 & 1.042869e+5 & 1.021905e+5 & 1.013128e+5 & 208.9 (531) \\\hline 

ALM    & 4.585705e+5  & 1.481379e+5 & 1.233182e+5 & 1.116683e+5 & 1.047410e+5 & 1.025611e+5 & 1.016589e+5 & 208.1 (581) \\\hline 

ISTA   & 2.345290e+5  & 1.267048e+5 & 1.137827e+5 & 1.079721e+5 & 1.040666e+5 & 1.025107e+5 & 1.018068e+5 & 196.8 (510) \\\hline 

SALSA  & 8.772957e+5  & 1.549462e+5 & 1.267379e+5 & 1.132676e+5 & 1.054600e+5 & 1.031346e+5 & 1.021898e+5 & 223.9 (663) \\\hline 

SADAL  & 2.524912e+5  & 1.271591e+5 & 1.133542e+5 & 1.068386e+5 & 1.021905e+5 & 1.004005e+5 & 9.961905e+4 & 113.5 (332) \\\hline
\end{tabular}
\end{center}}
\end{table}

\begin{table}[ht]{\scriptsize
\begin{center}\caption{Comparison of the algorithms for solving \eqref{prob:image-deconvolution} with $\rho=0.1$}\label{tab:comp-SALSA-3} \scriptsize
\begin{tabular}{|c|c|c|c|c|c|c|c|c|}\hline
solver & \multicolumn{7}{|c|}{obj in $k$-th iteration} & cpu (iter) \\\hline

       & 10 & 50 & 100 & 200 & 500 & 800 & 1000 & \\\hline

FALM-S & 9.868574e+5  & 8.771604e+5 & 8.487372e+5 & 8.271496e+5 & 8.110211e+5 & 8.065750e+5 & 8.050973e+5 & 37.7 (76) \\\hline 

FALM   & 9.876315e+5  & 8.629257e+5 & 8.369244e+5 & 8.210375e+5 & 8.097621e+5 & 8.067903e+5 & 8.058290e+5 & 25.7 (54) \\\hline 

FISTA  & 9.924884e+5  & 8.830263e+5 & 8.501727e+5 & 8.288459e+5 & 8.126598e+5 & 8.081259e+5 & 8.066060e+5 & 30.1 (79) \\\hline 

ALM-S  & 1.227588e+6  & 9.468694e+5 & 9.134766e+5 & 8.880703e+5 & 8.617264e+5 & 8.509737e+5 & 8.465260e+5 & 214.5 (537) \\\hline 

ALM    & 1.263787e+6  & 9.521381e+5 & 9.172737e+5 & 8.910902e+5 & 8.639917e+5 & 8.528932e+5 & 8.482666e+5 & 211.8 (588) \\\hline 

ISTA   & 1.048956e+6  & 9.396822e+5 & 9.161787e+5 & 8.951970e+5 & 8.700864e+5 & 8.589587e+5 & 8.541664e+5 & 293.8 (764) \\\hline 

SALSA  & 1.680608e+6  & 9.601661e+5 & 9.230268e+5 & 8.956607e+5 & 8.674579e+5 & 8.558580e+5 & 8.509770e+5 & 230.2 (671) \\\hline 

SADAL  & 1.060130e+6  & 9.231803e+5 & 8.956150e+5 & 8.735746e+5 & 8.509601e+5 & 8.420295e+5 & 8.383270e+5 & 112.5 (335) \\\hline
\end{tabular}
\end{center}}
\end{table}

\begin{figure}
\centering \subfigure{
\includegraphics[scale=0.5]{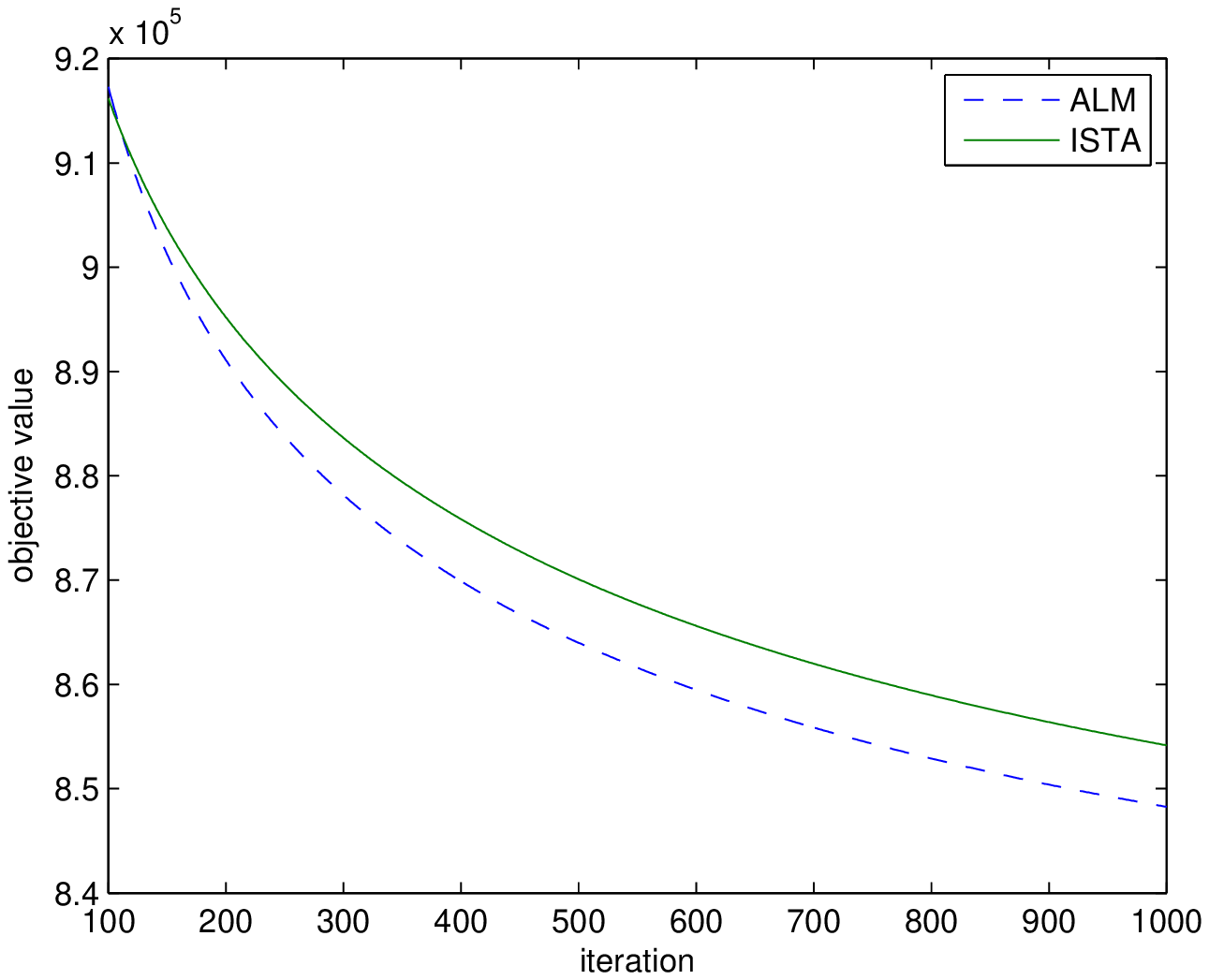}\label{fig:ALM-vs-ISTA}}
\centering \subfigure{
\includegraphics[scale=0.5]{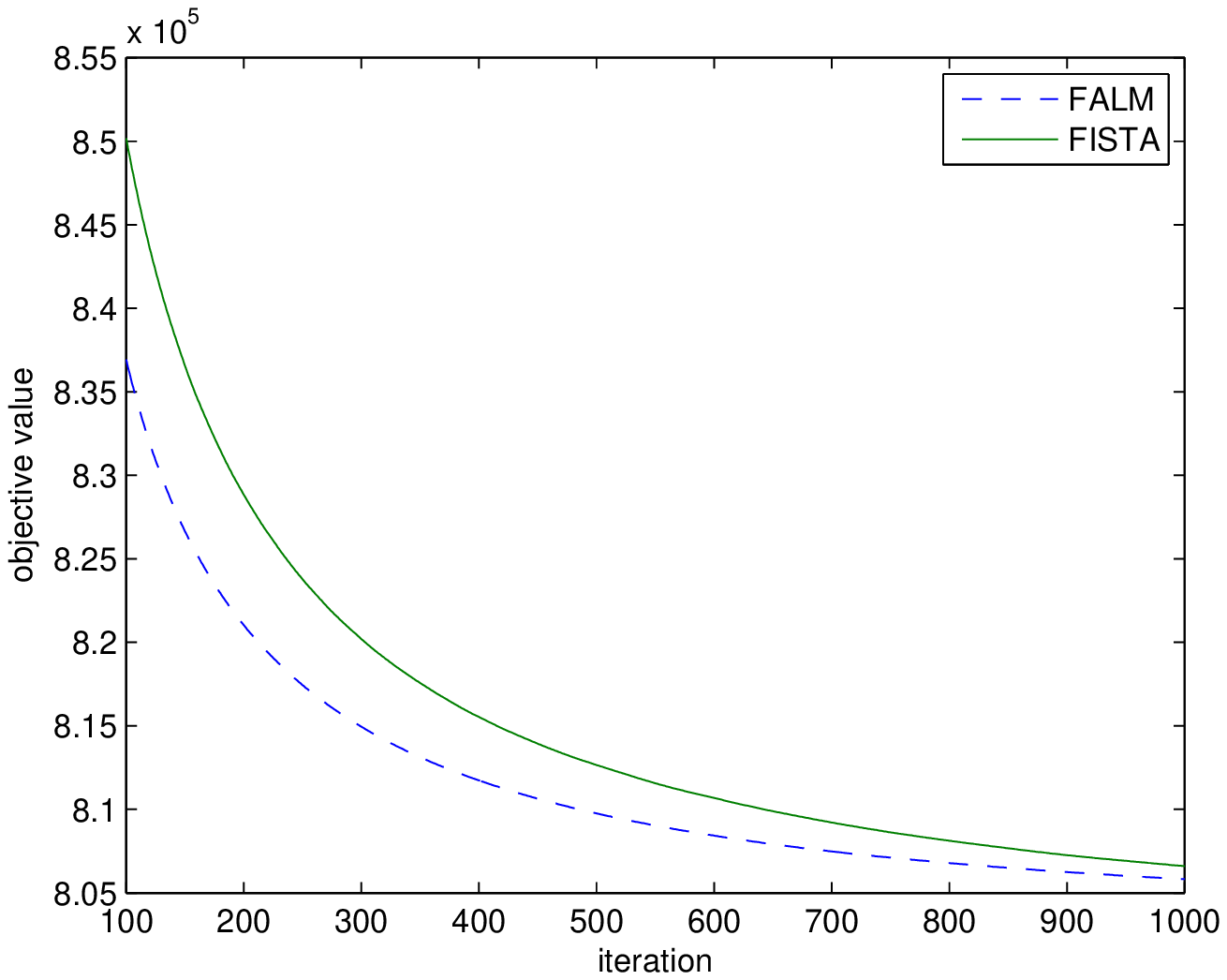}\label{fig:FALM_vs_FISTA}}
\centering \subfigure{
\includegraphics[scale=0.5]{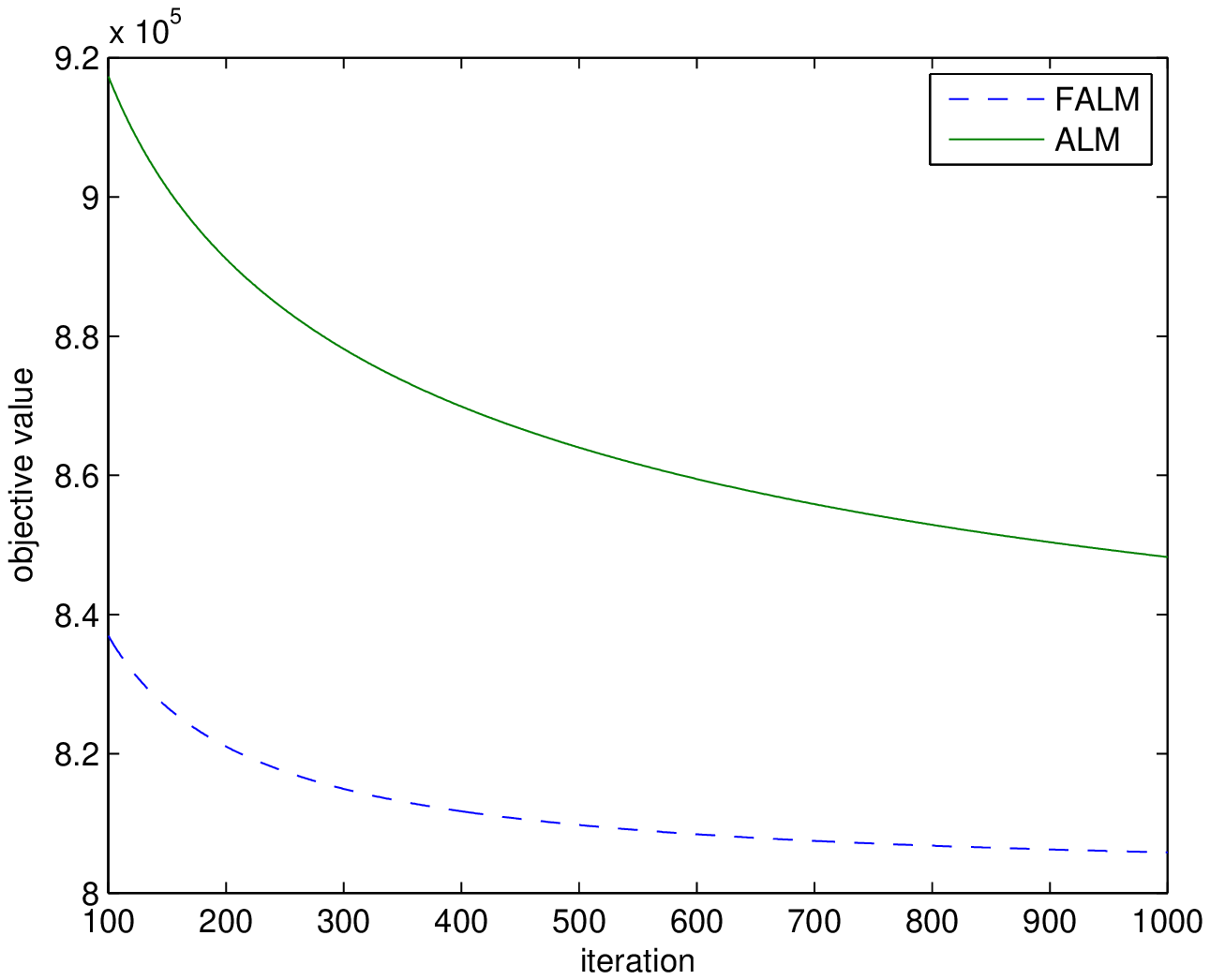}\label{fig:FALM_vs_ALM}}
\centering \subfigure{
\includegraphics[scale=0.5]{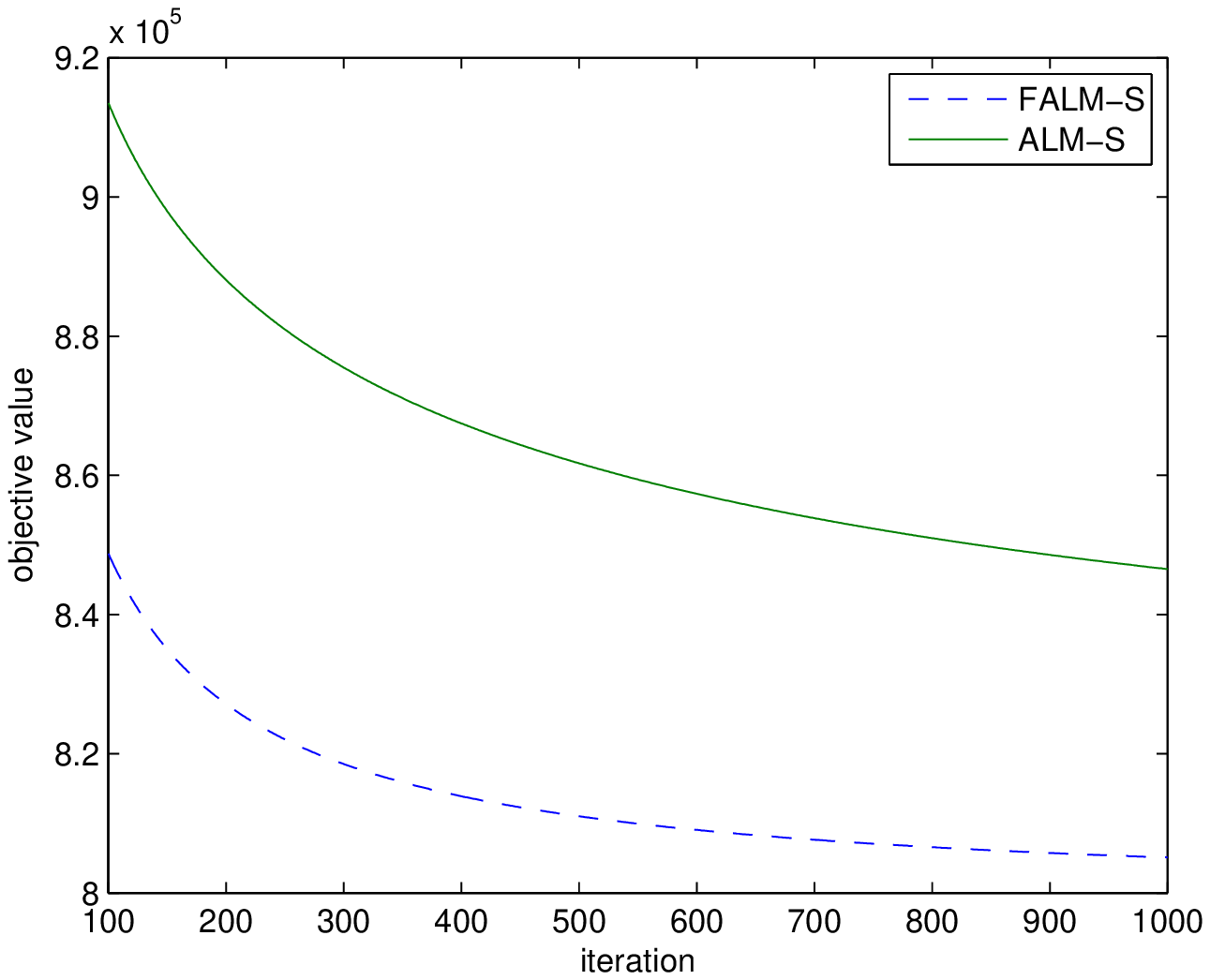}\label{fig:FALMs_vs_ALMs}}
\caption{comparison of the algorithms}\label{fig:comparison-all}
\end{figure}

From Tables \ref{tab:comp-SALSA-2} and \ref{tab:comp-SALSA-3} we see that in terms of the value of the objective function achieved after a specified number of iterations, the performance of FALM-S and FALM is always slightly better than that of FISTA and is much better than the performance of the other algorithms. On the two test problems, since FALM-S and FALM are always better than ALM-S and ALM, and FISTA is always better than ISTA, we can conclude that the Nesterov-type acceleration technique greatly speeds up the basic algorithms on these problems. Moreover, although sometimes in the early iterations FISTA (ISTA) is better than FALM-S and FALM (ALM-S and ALM), it is always worse than the latter two algorithms when the iteration number is large.
We also illustrate our comparisons graphically by plotting in Figure \ref{fig:comparison-all} the objective function value versus the number of iterations taken by these algorithms for solving \eqref{prob:image-deconvolution} with $\rho=0.1$. From Figure \ref{fig:comparison-all} we see clearly that for this problem, ALM outperforms ISTA, FALM outperforms FISTA,  FALM outperforms ALM and FALM-S outperforms ALM-S.

From the CPU times and the iteration numbers in the last columns of Tables \ref{tab:comp-SALSA-2} and \ref{tab:comp-SALSA-3} we see that, the fast versions are always much better than the basic versions of the algorithms. Since iterations of FISTA cost less than those of FALM-S (and FALM as well), we see that although FISTA takes 35\% (4\%) more iterations than FALM-S in the last column of Table \ref{tab:comp-SALSA-2} (\ref{tab:comp-SALSA-3}) it takes only 7\% more time (20\% less time).

We note that for the problems with $\rho = 0.01$ and $\rho = 0.1$, (i.e., for the results given in Tables \ref{tab:comp-SALSA-2} and \ref{tab:comp-SALSA-3}),
891 and 981, respectively, of the first 1000 iterations performed by FALM-S were skipping steps.
In contrast, none of the steps performed by ALM-S on either of these problems were skipping steps.
While the latter result is somewhat surprising, the fact that FALM-S performs many skipping steps is not, since the Nesterov-like acceleration approach is an {\it over-relaxation} approach that generates points that extrapolate beyond the previous point and the one produced by the ALM algorithm.


\section{Applications}\label{sec:application}
In this section, we describe how ALM and FALM can be applied to problems that can be formulated as RPCA problems to illustrate the use of Nesterov-type smoothing when the functions $f$ and $g$ do not satisfy the smoothness conditions required by the theorems in Sections \ref{sec:ALM} and \ref{sec:FALM}. Our numerical results show that our methods are able to solve huge problems that arise in practice; e.g., one problem involving roughly 40 million variables and 20 million linear constraints is solved in about three-quarters of an hour. We alse describe application of our methods to the SICS problem.

\subsection{Applications in Robust Principal Component Analysis}

In order to apply Algorithms \ref{alg:ALM} and \ref{alg:FALM} to \eqref{prob:RPCA}, we need to smooth both the nuclear norm $f(X):=\|X\|_*$ and the $\ell_1$ norm $g(Y):=\rho\|Y\|_1$. We again apply Nesterov's smoothing technique as in Section \ref{sec:comparison}. $g(Y)$ can be smoothed in the same way as the vector $\ell_1$ norm in Section \ref{sec:comparison}. We use $g_\sigma(Y)$ to denote the smoothed function with smoothness parameter $\sigma>0$. A smoothed approximation to $f(X)$ with smoothness parameter $\sigma>0$ is \bea\label{def:f-sigma-smooth} f_\sigma(X):=\max\{\langle X,W\rangle-\frac{\sigma}{2}\|W\|_F^2 : \|W\|\leq 1\}.\eea
It is easy to show that the optimal solution $W_\sigma(X)$ of \eqref{def:f-sigma-smooth} is \bea\label{def:W-sigma(X)}W_\sigma(X)=U\Diag(\min\{\gamma,1\})V^\top,\eea where $U\Diag(\gamma)V^\top$ is the singular value decomposition (SVD) of $X/\sigma.$ According to Theorem 1 in \citep{Nesterov-2005}, the gradient of $f_\sigma$ is given by $\nabla f_\sigma(X)=W_\sigma(X)$ and is Lipschitz continuous with Lipschitz constant $L(f_{\sigma})=1/\sigma$. After smoothing $f$ and $g$, we can apply Algorithms \ref{alg:ALM} and \ref{alg:FALM} to solve the following smoothed problem: \bea\label{prob:RPCA-Nuclear-L1-sigma}\min\{f_\sigma(X)+g_\sigma(Y):X+Y=M\}.\eea We have the following theorem about $\epsilon$-optimal solutions of problems \eqref{prob:RPCA} and \eqref{prob:RPCA-Nuclear-L1-sigma}.

\begin{theorem}\label{the:RPCA-Nuclear-L1-epsilon-optimal}
Let $\sigma=\frac{\epsilon}{2\max\{\min\{m,n\},mn\rho^2\}}$ and $\epsilon>0$. If $(X(\sigma),Y(\sigma))$ is an $\epsilon/2$-optimal solution to \eqref{prob:RPCA-Nuclear-L1-sigma}, then $(X(\sigma),Y(\sigma))$ is an $\epsilon$-optimal solution to \eqref{prob:RPCA}.
\end{theorem}

{\em Proof.} Let $D_f:=\max\{\half\|W\|_F^2 : \|W\|\leq 1\}=\half\min\{m,n\}$, $D_g:=\max\{\half\|Z\|_F^2 : \|Z\|_\infty\leq\rho\}=\half mn\rho^2$ and $(X^*,Y^*)$ and $(X^*(\sigma),Y^*(\sigma))$ be optimal solution to problems \eqref{prob:RPCA} and \eqref{prob:RPCA-Nuclear-L1-sigma}, respectively. Note that \bea\label{f-sigma-f-bound}f_\sigma(X)\leq f(X)\leq f_\sigma(X)+\sigma D_f, \forall X\in\br^\mtn\eea and \bea\label{g-sigma-g-bound}g_\sigma(Y)\leq g(Y)\leq g_\sigma(Y)+\sigma D_g,\forall Y\in\br^\mtn.\eea Using the inequalities in \eqref{f-sigma-f-bound} and \eqref{g-sigma-g-bound} and the facts that $(X(\sigma),Y(\sigma))$ is an $\epsilon/2$-optimal solution to \eqref{prob:RPCA-Nuclear-L1-sigma} and $\sigma\max\{D_f, D_g\}=\frac{\epsilon}{4}$, we have
\begin{align*} f(X(\sigma))+g(Y(\sigma))-f(X^*)-g(Y^*) & \leq  f_\sigma(X(\sigma))+g_\sigma(Y(\sigma))+\sigma D_f+\sigma D_g-f_\sigma(X^*)-g_\sigma(Y^*) \\ & \leq f_\sigma(X(\sigma))+g_\sigma(Y(\sigma))+\sigma D_f+\sigma D_g-f_\sigma(X^*(\sigma))-g_\sigma(Y^*(\sigma)) \\ & \leq \epsilon/2 + \sigma D_f+\sigma D_g \leq \epsilon/2 + \epsilon/4 + \epsilon/4 =  \epsilon. \qquad \endproof \end{align*}

Thus, according to Theorem \ref{the:RPCA-Nuclear-L1-epsilon-optimal}, to find an $\epsilon$-optimal solution to \eqref{prob:RPCA}, we need to find an $\epsilon/2$-optimal solution to \eqref{prob:RPCA-Nuclear-L1-sigma} with $\sigma=\frac{\epsilon}{2\max\{\min\{m,n\},mn\rho^2\}}$. We can either apply Algorithms \ref{alg:ALM} and \ref{alg:FALM} to solve \eqref{prob:RPCA-Nuclear-L1-sigma}, or apply Algorithms \ref{alg:ALM-S} and \ref{alg:FALM-S} to solve \eqref{prob:RPCA-Nuclear-L1-sigma} with only one functions (say $f(x)$) smoothed. The iteration complexity results in Theorems \ref{the:ALM-S} and \ref{the:FALM-S} hold since the gradients of $f_\sigma$ is Lipschitz continuous. However, the numbers of iterations needed by ALM and FALM to obtain an $\epsilon$-optimal solution to \eqref{prob:RPCA} become $O(1/\epsilon^2)$ and $O(1/\epsilon)$, respectively, due to the fact that the Lipschitz constant $L(f_\sigma)=1/\sigma= \frac{2\max\{\min\{m,n\},mn\rho^2\}}{\epsilon}=O(1/\epsilon)$.

The two subproblems at iteration $k$ of Algorithm \ref{alg:ALM} when applied to \eqref{prob:RPCA-Nuclear-L1-sigma} reduce to
\bea\label{sparse-low-rank-ALM-general-sub-1}X^{k+1}:=\arg\min_X f_\sigma(X)+g_\sigma(Y^k)+\langle\nabla g_\sigma(Y^k),M-X-Y^k\rangle+\oneotwomu\|X+Y^k-M\|_F^2,\eea
and
\bea\label{sparse-low-rank-ALM-general-sub-2}Y^{k+1}:=\arg\min_Y f_\sigma(X^{k+1})+\langle\nabla f_\sigma(X^{k+1}),M-X^{k+1}-Y\rangle+\oneotwomu\|X^{k+1}+Y-M\|_F^2+g_\sigma(Y).\eea
The first-order optimality conditions for \eqref{sparse-low-rank-ALM-general-sub-1} are: \bea\label{sparse-low-rank-ALM-general-sub-1-optcond}W_\sigma(X)-Z_\sigma(Y^k)+\oneomu(X+Y^k-M)=0,\eea where $W_\sigma(X)$ and $Z_\sigma(Y)$ are defined in \eqref{def:W-sigma(X)} and \eqref{def:Z-sigma(Y)}.
It is easy to check that
\begin{equation}\label{sparse-low-rank-ALM-general-sub-1-solution} X:=U\Diag(\gamma-\frac{\mu\gamma}{\max\{\gamma,\mu+\sigma\}})V^\top\end{equation} satisfies \eqref{sparse-low-rank-ALM-general-sub-1-optcond}, where
$U\Diag(\gamma)V^\top$ is the SVD of the matrix $\mu Z_\sigma(Y^k)-Y^k+M$. Thus, solving the subproblem \eqref{sparse-low-rank-ALM-general-sub-1} corresponds to an SVD. If we define $B:=\mu W_\sigma(X^{k+1})-X^{k+1}+M,$ it is easy to verify that \bea\label{sol:RPCA-Y}Y_{ij}=B_{ij}-\mu\min\{\rho,\max\{-\rho,\frac{B_{ij}}{\sigma+\mu}\}\}  \mbox{ for } i=1,\ldots,m \mbox{ and } j=1,\ldots,n\eea satisfies the first-order optimality conditions for \eqref{sparse-low-rank-ALM-general-sub-2}: $-W_\sigma(X^{k+1})+\oneomu(X^{k+1}+Y-M)+Z_\sigma(Y)=0.$ Thus, solving the subproblem \eqref{sparse-low-rank-ALM-general-sub-2} can be done very cheaply. The two subproblems at the $k$-th iteration of Algorithm \ref{alg:FALM} can be done in the same way and the main computational effort in each iteration of both ALM and FALM corresponds to an SVD.

\subsection{RPCA with Missing Data}
In some applications of RPCA, some of the entries of $M$ in \eqref{prob:RPCA} may be missing (e.g., in low-rank matrix completion problems where the matrix is corrupted by noise). Let $\Omega$ be the index set of the entries of $M$ that are observable and define the projection operator $\PCal_\Omega$ as: $(\PCal_\Omega(X))_{ij}=X_{ij}$, if $(i,j)\in\Omega$ and $(\PCal_\Omega(X))_{ij}=0$ otherwise.
It has been shown under some randomness hypotheses that the low rank $\bar{X}$ and sparse $\bar{Y}$ can be recovered with high probability by solving (see Theorem 1.2 in \citep{Candes-Li-Ma-Wright-RPCA-2009}),
\bea\label{prob:RPCA-Nuclear-L1-MC}(\bar{X},\bar{Y}):= \arg\min_{X,Y} \{ \|X\|_* + \rho\|Y\|_1 : \PCal_\Omega(X+Y) = \PCal_\Omega(M)\}. \eea
To solve \eqref{prob:RPCA-Nuclear-L1-MC} by ALM or FALM, we need to transform it into the form of \eqref{prob:RPCA}. For this we have
\begin{theorem}\label{the:RPCA-Nuclear-L1-MC-equiv}
$(\bar{X},\PCal_\Omega(\bar{Y}))$ is an optimal solution to \eqref{prob:RPCA-Nuclear-L1-MC} if
\bea\label{prob:RPCA-Nuclear-L1-MC-equiv}(\bar{X},\bar{Y})=\arg\min_{X,Y}\{ \|X\|_* + \rho\|\PCal_\Omega(Y)\|_1 : X + Y = \PCal_\Omega(M)\}. \eea
\end{theorem}
\begin{proof}
Suppose $(X^*,Y^*)$ is an optimal solution to \eqref{prob:RPCA-Nuclear-L1-MC}. We claim that $Y^*_{ij}=0,\forall (i,j)\notin\Omega$. Otherwise, $(X^*,\PCal_\Omega(Y^*))$ is feasible to \eqref{prob:RPCA-Nuclear-L1-MC} and has a strictly smaller objective function value than $(X^*,Y^*)$, which contradicts the optimality of $(X^*,Y^*)$. Thus, $\|\PCal_\Omega(Y^*)\|_1=\|Y^*\|_1$. Now suppose that $(\bar{X},\PCal_\Omega(\bar{Y}))$ is not optimal to \eqref{prob:RPCA-Nuclear-L1-MC}; then we have \bea\label{proof:the:RPCA-Nuclear-L1-MC-equiv-eq1}\|X^*\|_*+\rho\|\PCal_\Omega(Y^*)\|_1 = \|X^*\|_*+\rho\|Y^*\|_1 < \|\bar{X}\|_* +\rho\|\PCal_\Omega(\bar{Y})\|_1. \eea By defining a new matrix $\tilde{Y}$ as \beaa\tilde{Y}_{ij}=\left\{\ba{ll} Y^*_{ij}, & (i,j)\in\Omega \\ -X^*_{ij}, & (i,j)\notin\Omega, \ea\right.\eeaa we have that $(X^*,\tilde{Y})$ is feasible to \eqref{prob:RPCA-Nuclear-L1-MC-equiv} and $\|\PCal_\Omega(\tilde{Y})\|_1=\|\PCal_\Omega(Y^*)\|_1$. Combining this with \eqref{proof:the:RPCA-Nuclear-L1-MC-equiv-eq1}, we obtain \beaa\|X^*\|_*+\rho\|\PCal_\Omega(\tilde{Y})\|_1<\|\bar{X}\|_*+\rho\|\PCal_\Omega(\bar{Y})\|_1,\eeaa which contradicts the optimality of $(\bar{X},\bar{Y})$ to \eqref{prob:RPCA-Nuclear-L1-MC-equiv}. Therefore, $(\bar{X},\PCal_\Omega(\bar{Y}))$ is optimal to \eqref{prob:RPCA-Nuclear-L1-MC}.
\end{proof}

The only differences between \eqref{prob:RPCA} and \eqref{prob:RPCA-Nuclear-L1-MC-equiv} lie in that the matrix $M$ is replaced by $\PCal_\Omega(M)$ and $g(Y)=\rho\|Y\|_1$ is replaced by $\rho\|\PCal_\Omega(Y)\|_1$. A smoothed approximation $g_\sigma(Y)$ to $g(Y):=\rho\|\PCal_\Omega(Y)\|_1$ is given by
\bea\label{smooth-L1-P_Omega(Y)}g_\sigma(Y):=\max\{\langle \PCal_\Omega(Y),Z\rangle - \frac{\sigma}{2}\|Z\|_F^2: \|Z\|_\infty\leq\rho\},\eea and \bea\label{smooth-L1-P_Omega(Y)-optsol}(\nabla g_\sigma(Y))_{ij}= \min\{\rho,\max\{(\PCal_\Omega(Y))_{ij}/\sigma,-\rho\}\},\mbox{ for } 1\leq i\leq m \mbox{ and } 1\leq j\leq n.\eea According to Theorem 1 in \citep{Nesterov-2005}, $\nabla g_\sigma(Y)$ is Lipschitz continuous with $L_\sigma(g)=1/\sigma$. Thus the convergence and iteration complexity results in Theorems \ref{the:ALM-S} and \ref{the:FALM-S} apply.
The only changes in Algorithms \ref{alg:ALM} and \ref{alg:ALM-S} and Algorithms \ref{alg:FALM} and \ref{alg:FALM-S} are: replacing $M$ by $\PCal_\Omega(M)$ and computing $Y^{k+1}$ using \eqref{sol:RPCA-Y} with $B$ is replaced by $\mathcal{P}_\Omega(B)$.

\subsection{Numerical Results on RPCA Problems} In this section, we report numerical results obtained using the ALM method to solve RPCA problems with both complete and incomplete data matrices $M$. We compare the performance of ALM with the exact ADM (EADM) and the inexact ADM (IADM) methods in \citep{Lin-Chen-Wu-Ma-ADM-RPCA-2009}. The MATLAB codes of EADM and IADM were downloaded from $http://watt.csl.illinois.edu/\sim perceive/matrix-rank/sample\_code.html$ and their default settings were used. To further accelerate ALM, we adopted the continuation strategy used in EADM and IADM. Specifically, we set $\mu_{k+1}:=\max\{\bar{\mu},\eta\mu_k\}$, where $\mu_0=\|M\|/1.25, \bar{\mu}=10^{-6}$ and $\eta=2/3$ in our numerical experiments. Although in some iterations this violates the requirement $\mu \leq \min\{\frac{1}{L(f_\sigma)},\frac{1}{L(g_\sigma)}\}$ in Corollaries \ref{cor:ALM-nonsmooth} and \ref{cor:FALM-nonsmooth}, we see from our numerical results reported below that ALM and FALM still work very well. We also found that by adopting this updating rule for $\mu$, there was not much difference between the performance of ALM and that of FALM. So we only compare ALM with EADM and IADM. As in Section \ref{sec:comparison}, since we applied ALM to a smoothed problem, we implemented ALM as SADAL. The initial point in ALM was set to $(X^0,Y^0) = (M,\mathbf{0})$ and the initial Lagrange multiplier was set to $\Lambda^0=-\nabla g_\sigma(Y^0)$. We set the smoothness parameter $\sigma=10^{-6}$. Solving subproblem \eqref{sparse-low-rank-ALM-general-sub-1} requires computing an SVD (see \eqref{sparse-low-rank-ALM-general-sub-1-solution}). However, we do not have to compute the whole SVD, as only the singular values that are larger than the threshold $\tau=\frac{\mu\gamma}{\max\{\gamma,\mu+\sigma\}}$ and the corresponding singular vectors are needed. We therefore use PROPACK \citep{Larsen-Propack}, which is also used in EADM and IADM, to compute these singular values and corresponding singular vectors. To use PROPACK, one has to specify the number of leading singular values (denoted by $sv_k$) to be computed at iteration $k$. We here adopt the strategy suggested in \citep{Lin-Chen-Wu-Ma-ADM-RPCA-2009} for EADM and IADM. This strategy starts with $sv_0=100$ and updates $sv_k$ via:
\begin{eqnarray*}
sv_{k+1} = \left\{
\begin{array}{ll}
svp_k + 1, & \mbox{ if } svp_k < sv_k \\
\min\{svp_k+round(0.05d),d\}, & \mbox{ if } svp_k = sv_k,
\end{array}\right.
\end{eqnarray*}
where $d=\min\{m,n\}$ and $svp_k$ is the number of singular values that are larger than the threshold $\tau$.

In all our experiments $\rho$ was chosen equal to $1/\sqrt{m}$. We stopped ALM, EADM and IADM when the relative infeasibility was less than $10^{-7}$, i.e., $\|X+Y-M\|_F < 10^{-7}\|M\|_F$.

\subsubsection{Background Extraction from Surveillance Video}
Extracting the almost still background from a sequence frames of video is a basic task in video surveillance. This problem is difficult due to the presence of
moving foregrounds in the video. Interestingly, as shown in \citep{Candes-Li-Ma-Wright-RPCA-2009}, this problem can be formulated as a RPCA problem \eqref{prob:RPCA}. By stacking the columns of each frame into a long vector, we get a matrix $M$ whose columns correspond to the sequence of frames of the video.
This matrix $M$ can be decomposed into the sum of two matrices $M:=\bar{X}+\bar{Y}$. The matrix $\bar{X}$, which represents the background in the frames, should be of low rank due to the correlation between frames. The matrix $\bar{Y}$, which represents the moving objects in the foreground in the frames, should be sparse since
these objects usually occupy a small portion of each frame. We apply ALM to solve \eqref{prob:RPCA} for two videos introduced in \citep{Li-Huang-Gu-Tian-2004}.

\begin{figure}\hspace{-4cm}
\centering \subfigure{
\includegraphics[scale=0.5]{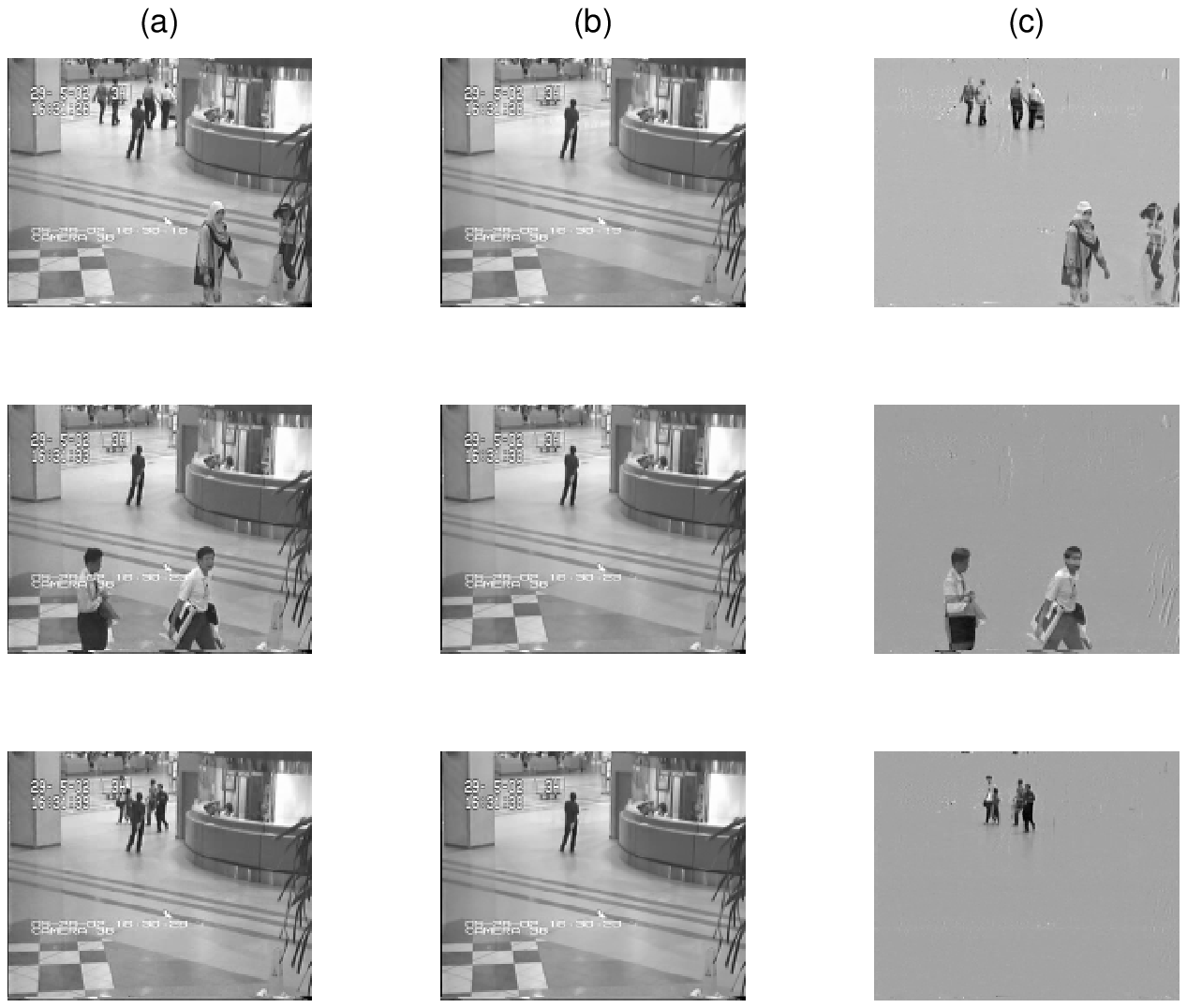}\label{fig:Hall_airport}}\hspace{-1cm}
\centering \subfigure{
\includegraphics[scale=0.5]{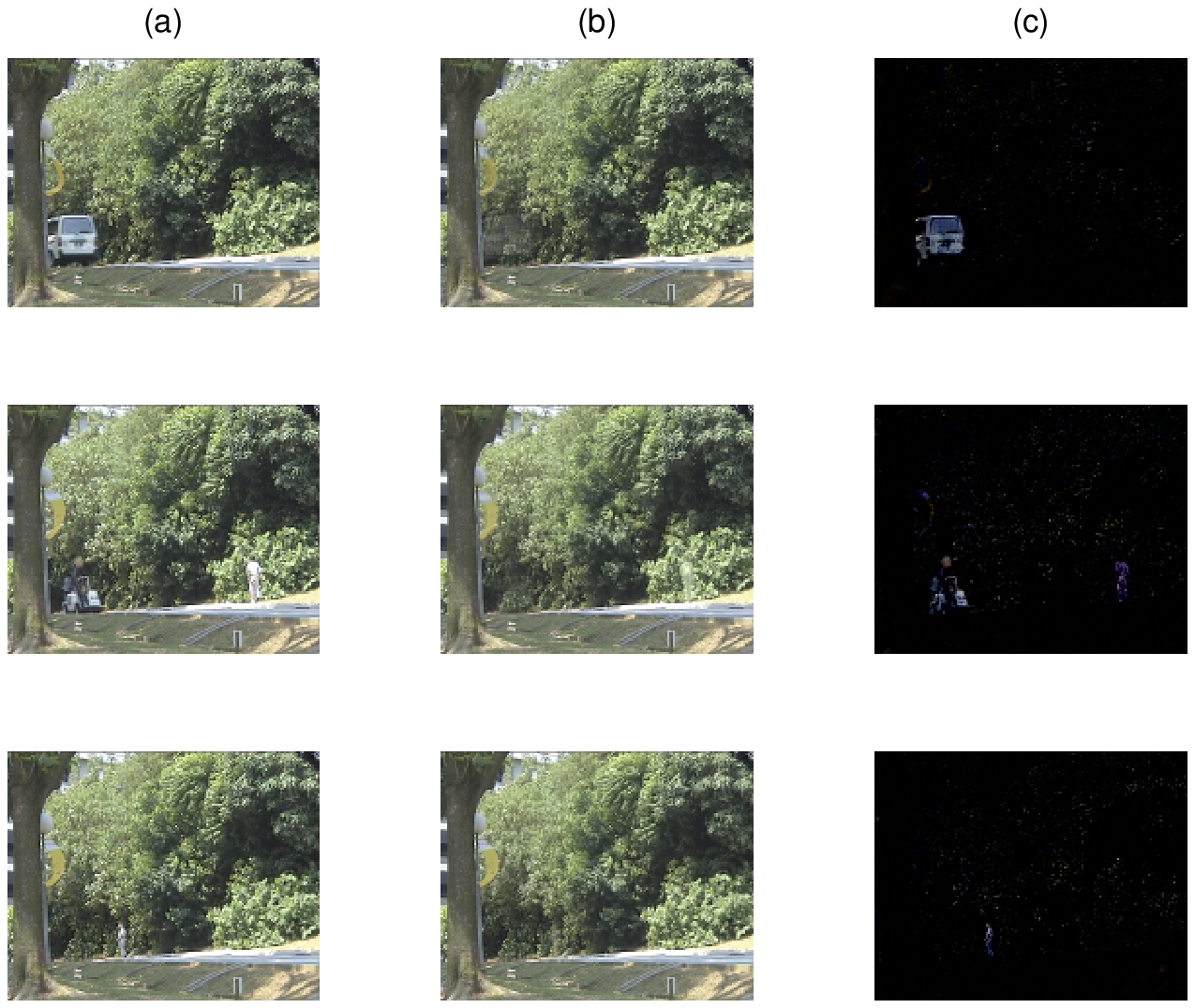}\label{fig:Escalator_airport}}\hspace{-4cm}
\label{fig:video}\caption{In the first 3 columns: (a) Video sequence. (b) Static background recovered by our ALM. Note that the man who kept still in the 200 frames stays as in the background. (c) Moving foreground recovered by our ALM. In the last 3 columns: (a) Video sequence. (b) Static background recovered by our ALM. (c) Moving foreground recovered by our ALM.}
\end{figure}

Our first example is a sequence of 200 grayscale frames of size $144\times 176$ from a video of a hall at an airport. Thus the matrix $M$ is in $\br^{25344\times 200}$. The second example is a sequence of 320 color frames from a video taken at a campus. Since the video is colored, each frame is an image stored in the RGB format, which is a $128\times 160 \times 3$ cube. The video is then reshaped into a $128\times 160$ by $3\times 320$ matrix, i.e., $M\in\br^{20480\times 960}$. Some frames of the videos and the recovered backgrounds and foregrounds are shown in Figure \ref{fig:video}. We only show the frames produced by ALM, because EADM and IADM produce visually identical results. From these figures we can see that ALM can effectively separate the nearly still background from the moving foreground. Table \ref{tab:surveillance-video} summarizes the numerical results on these problems. The CPU times are reported in the form of $hh:mm:ss$.
From Table \ref{tab:surveillance-video} we see that although ALM is slightly worse than IADM, it is much faster than EADM in terms of both the number of SVDs and CPU times. We note that the numerical results in \citep{Candes-Li-Ma-Wright-RPCA-2009} show that the model \eqref{prob:RPCA} produces much better results than other competing models for background extraction in surveillance video.

\begin{table}[ht]{\scriptsize
\begin{center}\caption{Comparison of ALM and EADM on surveillance video problems}\label{tab:surveillance-video} \small
\begin{tabular}{|l c c | c r| c r| c r |}\hline
\multicolumn{3}{|c|}{} &
\multicolumn{2}{|c|}{Exact ADM} & \multicolumn{2}{|c|}{Inexact ADM} & \multicolumn{2}{|c|}{ALM}  \\\hline

Problem & $m$ & $n$ & SVDs & CPU & SVDs & CPU & SVDs & CPU \\\hline

Hall (gray)     & 25344 & 200  & 550 & 40:15  & 38 & 03:47 &  43 & 04:03  \\\hline

Campus (color) & 20480 & 960 &  651 & 13:54:38 & 40 & 43:35  & 46 & 46:49  \\\hline
\end{tabular}
\end{center}}
\end{table}

\subsubsection{Random Matrix Completion Problems with Grossly Corrupted Data}
For the matrix completion problem \eqref{prob:RPCA-Nuclear-L1-MC}, we set $M:=A+E$, where the rank $r$ matrix $A\in\br^{n\times n}$ was created as the product $A_LA_R^\top$, of random matrices $A_L\in\br^{n\times r}$ and $A_R\in\br^{n\times r}$ with i.i.d. Gaussian entries $\mathcal{N}(0,1)$ and the sparse matrix $E$ was generated by choosing its support uniformly at random and its nonzero entries uniformly i.i.d. in the interval $[-500,500]$. In Table \ref{tab:MC}, $rr:=\rank(A)/n$, $spr:=\|E\|_0/n^2$, the relative errors $relX := \|X-A\|_F / \|A\|_F$ and $relY: = \|Y-E\|_F / \|E\|_F$, and the sampling ratio of $\Omega$, $SR=m/n^2$. The $m$ indices in $\Omega$ were generated uniformly at random. We set $\rho=1/\sqrt{n}$ and stopped ALM when the relative infeasibility $\|X+Y-\PCal_\Omega(M)\|_F/\|\PCal_\Omega(M)\|_F < 10^{-5}$ and for our continuation strategy, we set $\mu_0 = \|\PCal_\Omega(M)\|_F/1.25$. The test results obtained using ALM to solve \eqref{prob:RPCA-Nuclear-L1-MC-equiv} with the nonsmooth functions replaced by their smoothed approximations are given in Table \ref{tab:MC}. From Table \ref{tab:MC} we see that ALM recovered the test matrices from a limited number of observations. Note that a fairly high number of samples was needed to obtain small relative errors due to the presence of noise. The number of iterations needed was almost constant (around 36), no matter the size of the problems. The CPU times (in seconds) needed are also reported.
\begin{table}[ht]{
\begin{center}\caption{{Numerical results for noisy matrix completion problems}}\label{tab:MC}
\begin{tabular}{|c c| c c c c |c c c c |}\hline
$rr$ & $spr$    & iter & relX & relY & cpu & iter & relX & relY & cpu            \\\hline
\multicolumn{2}{|c|}{} & \multicolumn{4}{|c|}{$SR=90\%, n = 500$} & \multicolumn{4}{|c|}{$SR=80\%, n=500$} \\\hline

$0.05$ & $0.05$ &  36 & 4.60e-5 & 4.25e-6 & 137 & 36  & 3.24e-5 & 4.31e-6 & 153  \\\hline
$0.05$ & $0.1$  &  36 & 4.68e-5 & 5.29e-6 & 156 & 36  & 4.40e-5 & 4.91e-6 & 161  \\\hline
$0.1$  & $0.05$ &  36 & 4.04e-5 & 3.74e-6 & 128 & 36  & 1.28e-3 & 1.33e-4 & 129  \\\hline
$0.1$  & $0.1$  &  36 & 6.00e-4 & 4.50e-5 & 129 & 35  & 1.06e-2 & 7.59e-4 & 124  \\\hline
\multicolumn{2}{|c|}{} & \multicolumn{4}{|c|}{$SR=90\%, n = 1000$} & \multicolumn{4}{|c|}{$SR=80\%, n=1000$} \\\hline

$0.05$ & $0.05$ & 37 & 3.10e-5 &  3.96e-6 & 1089 & 37 & 2.27e-5 & 4.14e-6 & 1191 \\\hline
$0.05$ & $0.1$  & 37 & 3.20e-5 &  4.93e-6 & 1213 & 37 & 3.00e-5 & 4.66e-6 & 1271 \\\hline
$0.1$  & $0.05$ & 37 & 2.68e-5 &  3.34e-6 & 982  & 37 & 1.75e-4 & 2.49e-5 & 994 \\\hline
$0.1$  & $0.1$  & 37 & 3.64e-5 &  4.51e-6 & 1004 & 36 & 4.62e-3 & 4.63e-4 & 965 \\\hline
\end{tabular}
\end{center}}
\end{table}

\subsection{Sparse Inverse Covariance Selection}

In  \citep{Scheinberg-Ma-Goldfarb-NIPS-2010} ALM method was successfully
applied to the Sparse Inverse Covariance Selection problem:
\bea\label{prob:min-f-g-SICS} \min_{X\in S^n_{++}} \quad F(X) \equiv f(X) + g(X), \eea where $f(X)=-\log\det(X)+\langle S,X\rangle$ and $g(X)=\rho\|X\|_1$.

 Note that in our case $f(X)$ does not have Lipschitz continuous
gradient in general. Moreover, $f(X)$ is
 only defined for positive definite matrices while
$g(X)$  is defined everywhere. These properties of the  objective function  make the
SICS problem especially challenging for optimization methods.  Nevertheless, we can still
 apply  Algorithm \ref{alg:ALM-S} and obtain the complexity bound in  Theorem \ref{the:ALM-S}  as follows.
As proved in \citep{Lu-covsel-siopt-2009}, the optimal solution $X^*$ of \eqref{prob:min-f-g-SICS} satisfies $X\succeq \alpha I$, where $\alpha = \frac{1}{\|S\|+n\rho},$
(see Proposition 3.1 in \citep{Lu-covsel-siopt-2009}).
 Therefore, the SICS problem  \eqref{prob:min-f-g-SICS} can be formulated as:
\bea\label{prob:min-f-g-constraint}\min_{X,Y} \{ f(X) + g(Y) : X - Y = 0, X \in \mathcal{C}, Y \in \mathcal{C}\},\eea where $\mathcal{C}:=\{X\in S^n: X\succeq\frac{\alpha}{2} I\}$. We can apply   Algorithm \ref{alg:ALM-S} and
  Theorem \ref{the:ALM-S} as per Remark \ref{rem:convex-set}. The difficulty
arises, however, when performing minimization in $Y$ (Step 5 of  Algorithm \ref{alg:ALM-S}) with the constraint  $Y \in \mathcal{C}$. Without this constraint,
the minimization is obtained by a matrix shrinkage operation. However, the problem becomes harder to solve with this additional constraint.
Minimization in $X$ (Step 3 of  Algorithm \ref{alg:ALM-S}) with or without
the constraint  $X \in \mathcal{C}$ is accomplished by performing an SVD of the current iterate $Y^k$. Hence the constraint can be easily imposed. Also note that once the SVD is computed both $\nabla f(X^{k+1})$ and $\nabla f(Y^k)$ are readily available (see \citep{Scheinberg-Ma-Goldfarb-NIPS-2010} for details). This implies that either skipping or nonskipping iterations of  Algorithm \ref{alg:ALM-S} can be performed at the same cost as one ISTA iteration.

Instead of imposing constraint  $Y \in \mathcal{C}$ in Step 5 of  Algorithm \ref{alg:ALM-S} we can obtain feasible solutions by a line search on $\mu$.
 We know that the constraint $ X\succeq\frac{\alpha}{2} I$ is not tight at the solution. Hence if we start
the algorithm with $X\succeq\alpha I $ and restrict the step size $\mu$ to be sufficiently small
then the iterates of the method will remain  in $\mathcal{C}$.
Similarly, one can apply ISTA with small steps to remain in $\mathcal{C}$.
Note however, that the bound on the Lipschitz constant of the gradient of $f(X)$ is $1/\alpha^2$ and hence can be very large. It is not practical to restrict $\mu$ in the algorithm to be smaller than $\alpha^2$, since $\mu$ determines the step size at each iteration.
The advantage of ALM methods  over ISTA in this case is that as soon as the  $Y \in \mathcal{C}$ is relaxed ISTA can no longer be applied,
while  ALM/SADAL  can be applied  and indeed works very well. The theory
in this case only applies once certain proximity to the optimal solution
has been reached. But as shown in \citep{Scheinberg-Ma-Goldfarb-NIPS-2010}, the SADAL method is computationally superior to other state-of-the-art methods for SICS.

We have also applied the FALM method to the SICS problem, but we have not observed any advantage  over ALM for this particular application.

\section{Conclusion}\label{sec:conclude}

In this paper, we proposed both basic and accelerated versions of alternating linearization methods for minimizing the sum of two convex functions. Our basic methods require at most $O(1/\epsilon)$ iterations to obtain an $\epsilon$-optimal solution, while our accelerated methods require at most $O(1/\sqrt{\epsilon})$ iterations with only a small additional amount of computational effort at each iteration. Numerical results on image deblurring, background extraction from surveillance video and matrix completion with grossly corrupted data are reported. These results demonstrate the efficiency and the practical potential of our algorithms.

\section*{Acknowledgement}
We would like to thank Dr. Zaiwen Wen for insightful discussions on the topic of this paper.

\bibliographystyle{siam}
\bibliography{C:/Mywork/Optimization/work/reports/bibfiles/All}
\end{document}